\documentclass[12pt]{amsart}
\usepackage{a4wide}
\usepackage[utf8]{inputenc}
\usepackage{amsmath}
\usepackage{amssymb}
\usepackage{amsfonts}
\usepackage{amsopn}
\usepackage{graphicx}
\usepackage{enumerate}
\usepackage[dvipsnames]{xcolor}
\usepackage{mathtools}
\usepackage{tikz-cd}
\usepackage{mathrsfs}
\usepackage{bbm}
\usepackage{yhmath}
\usepackage{cite}
\usepackage[colorlinks, linkcolor=blue, citecolor=blue, pagebackref]{hyperref}

\newtheorem{theorem}{Theorem}[section]
\newtheorem{lemm}[theorem]{Lemma}
\newtheorem{prop}[theorem]{Proposition}
\newtheorem{coro}[theorem]{Corollary}

\theoremstyle{definition}
\newtheorem{defi}[theorem]{Definition}

\newtheorem{conj}[theorem]{Conjecture}
\theoremstyle{remark}
\newtheorem{remark}[theorem]{Remark}
\numberwithin{equation}{section}
\allowdisplaybreaks

\numberwithin{equation}{section}

\begin{document}
\title[Two-parameter quantum groups vs One-parameter quantum groups]
{Two-parameter versus one-parameter\\
quantum groups: structural differences\\
arising from the type~$D$ phenomenon}

\author[Hu]{Naihong Hu}
\address{School of Mathematical Sciences, MOE-KLMEA \& SH-KLPMMP, East China Normal University,
MinHang Campus, DongChuan Road 500, Shanghai 200241, PR China}
\email{nhhu@math.ecnu.edu.cn}

\author[Rosso]{Marc Rosso}
\address{UFR Mathematics de Université Paris Cité, CNRS, IMJ-PRG, Bâtiment Sophie Germain – 8 place Aurélie
Nemours, 75013 Paris}
\email{marc.rosso@imj-prg.fr}

\author[Xu]{Xiao Xu{$^*$}}
\address{School of Mathematical Sciences, MOE-KLMEA \& SH-KLPMMP, East China Normal University,
MinHang Campus, DongChuan Road 500, Shanghai 200241, PR China}
\email{52275500004@stu.ecnu.edu.cn}

\thanks{$^*$Corresponding author. This work is
	supported by the NNSF of China (Grant No. 12171155), and in part by the Science and Technology Commission of Shanghai Municipality (Grant No. 22DZ2229014).}



\begin{abstract}
This paper investigates several differences between two-parameter and one-parameter quantum (affine) algebras.
We compare their Hopf algebra structures, algebra structures, automorphism groups of positive parts and extended Borel algebras,
and module categories beyond the finite-dimensional setting.
These differences are closely related to what we call \textit{``the type $D$ phenomenon''},
which widely appears in the structure of both finite and affine two-parameter quantum groups.
\end{abstract}

\maketitle
\tableofcontents

\section{Introduction}
\textbf{1.1}
Two-parameter quantum groups were first introduced by Takeuchi \cite{Takeuchi 1990} as a natural generalization of the Drinfeld-Jimbo
quantum groups.
Motivated by two-parameter deformations arising from down-up algebras \cite{BenkartW2001},
Benkart-Witherspoon and their collaborators systematically investigated the structure and representation theory of two-parameter quantum groups of type $A$ \cite{BW FIC 2004,BW2,BenkartW2004,Benkart06}.
Subsequently, Bergeron-Gao-Hu extended the theory to the remaining classical types and established their Drinfeld double structures together with the corresponding finite-dimensional weight representation theories \cite{BGH 2006,BGH 2007}.
For exceptional types, analogous constructions were obtained by Bai-Hu, Chen-Hu-Wang and Hu-Shi \cite{BaiH2008,CHW F4,HS G2}.

A fundamental question is to what extent two-parameter quantum groups differ from their one-parameter counterparts.
From a deformation-theoretic viewpoint,
Hu-Pei \cite{HuPei2008,HuPeiCA 2012} showed that two-parameter quantum groups arise from suitable $2$-cocycle
deformations of one-parameter quantum groups.
Moreover, the categories of finite-dimensional weight modules are equivalent as braided tensor categories.

Nevertheless, a number of essential differences have been discovered.
For instance, Lusztig symmetries fail to exist in general when $\mathrm{rank}\,\mathfrak{g}>2$ \cite{BGH 2006},
indicating that two-parameter quantum groups possess substantially less symmetry than their one-parameter counterparts.
At roots of unity, suitable choices of $r=q^x$ and $s=q^y$, where $q$ is a primitive $\ell$-th root of unity and $x,y\in\mathbb{Z}$,
give rise to exotic small quantum groups that are not isomorphic, as Hopf algebras, to the standard one-parameter small quantum groups.
Even when restricted to $\ell\in\{4,5,6,7,8\}$,
the exotic small quantum groups counted in \cite{HX} comprise approximately five times as many Hopf algebra isomorphism classes as their standard counterparts.
More recently, Hu-Wang \cite{HuWang 2026} established the Harish-Chandra isomorphism theorem for two-parameter quantum groups when
$\mathfrak{g}$ is simple of even rank,
and identified new central generators in the odd rank case.

The distinctions between two-parameter and one-parameter quantum groups also extend to the affine setting.
Hu-Rosso-Zhang \cite{HRZ 2008} constructed the quantum affine Lyndon basis of $U_{r,s}(\widehat{\mathfrak{sl}_n})$ and established the Drinfeld realization through a combinatorial approach which differs substantially from Beck's braid-group method in the one-parameter setting.
Subsequent developments include generating functions with $\tau$-invariance, vertex representations, Drinfeld realizations of
twisted affine algebras, and $R$-matrix realizations of affine
and super quantum algebras \cite{HuZhang2014,HuZhang JA,JingLJA 2017,JingZ09,JingZ JMP2016,
HXZ,ZhHJ C,ZHX D,ZhangH 2018,HJZ JMP2024}.
These results further demonstrate that the two-parameter theory possesses structural features that cannot always be reduced to the
one-parameter framework.

\textbf{1.2}
The main goal of this paper is to clarify the differences between two-parameter and one-parameter quantum algebras,
not only at the level of Hopf structures, but also, more fundamentally,
in their algebra structures and the corresponding representation categories.

Benkart-Whiterspoon studied the Hopf algebraic distinctions in type $A$, see \cite{BenkartW2004}.
They proved that for $n>2$, $U_{r, s}(\mathfrak{sl}_{n+1})\cong U_{q, q^{-1}}(\mathfrak{sl}_{n+1})$ as Hopf algebras if and only if $(r, s)=(q, q^{-1})$.
However, they overlooked the case $(r, s)=(q^{-1}, q)$, which also yields an isomorphism between these two Hopf algebras.
We fill this gap and further generalize the results for type $A$ to all remaining finite types as well as affine types, see Theorem \ref{thm isofin} and Corollary \ref{cor iso}.

More importantly, the differences are also reflected in the automorphism groups of their positive parts.
Andruskiewitsch-Dumas \cite{Andrus08} conjectured that the automorphism group of $U_q^+(\mathfrak g)$ is generated by torus automorphisms and those induced by Dynkin diagram automorphisms. They verified this conjecture in the case of type $A_2$, and it was later proved by Yakimov \cite{Yakimov 2014}.
However, the absence of symmetry makes two-parameter quantum groups structurally more rigid.
Using $2$-cocycle deformation theory \cite{HuPei2008} together with Yakimov's results,
we give a characterization of the automorphism group of $U_{r,s}^+(\mathfrak{g})$.
By comparing these groups, we conclude that if $\mathfrak{g} = \mathfrak{sl}_{n+1}$, $\mathfrak{so}_{2n}$ or $E_6$, then $U_{r,s}^+(\mathfrak{g}) \not\cong U_q^+(\mathfrak{g})$ as algebras.
Furthermore, to compensate for the broken symmetry in the two-parameter setting,
we introduce the notion of extended Borel algebras, see Definition \ref{def extended}.
In this framework, automorphisms arising from Dynkin diagram symmetries are restored, as in the one-parameter case,
see Propositions \ref{aut extendedAn} and \ref{aut extendedD4}.
In type $D_4$, we prove that some diagram automorphisms cannot be realized at the level of the root lattice, see Proposition \ref{prop aut D4}.
To overcome this obstruction, we pass from the root lattice to the weight lattice and consider the corresponding weight-type extended Borel algebra.
We show that the full Dynkin diagram symmetry is recovered if and only if the extended Borel algebra is of weight-type, see Theorem \ref{prop aut D4weight}.

One may naturally ask how the additional parameter affects the structure
of quantum groups.
The factor $rs$, which becomes $1$ under the one-parameter specialization
$r=q$, $s=q^{-1}$, occurs in the structure constants and commutation
relations of two-parameter quantum groups.
For example, in $U_{r,s}(\mathfrak{so}_{2n})$, we have
$e_{n-1}e_n=rs\,e_ne_{n-1}$,
whereas the corresponding generators commute in
$U_{q,q^{-1}}(\mathfrak{so}_{2n})$.
We refer to this change from quasi-commutativity to commutativity
under the one-parameter specialization as the \textit{type $D$ phenomenon}
(cf. \cite{BGH 2006}).
Similar behavior occurs between composite quantum root vectors
in other finite types and also appears in a broad range of
two-parameter quantum affine algebras
(cf. \cite{HRZ 2008, HuZhang2014, Zhang phd 2007}).

Table~\ref{tab type D examples} lists examples in finite and affine types.
Every displayed relation becomes a commutative relation under
the specialization $r=q$, $s=q^{-1}$.
The references in the last column specify the notation and conventions;
additional choices of root vectors are given below the table.

\begin{table}[htbp]
  \centering
  \caption{Examples of the type $D$ phenomenon in finite and affine types.}
  \label{tab type D examples}
  \begingroup
  \small
  \renewcommand{\arraystretch}{1.7}
  \setlength{\tabcolsep}{5pt}
  \begin{tabular}{lll}
    \hline
    Type & Commutation relation & Notation and conventions \\
    \hline

    $A_n$, $n\geq 3$
    & $\mathcal{E}_{2,1}Y=rs\,Y\mathcal{E}_{2,1}$
    & \cite{BW FIC 2004} \\[3pt]

    $B_n$, $n\geq 3$
    & $E_{n-2,(n-1)'}e_n
       =(rs)^2e_nE_{n-2,(n-1)'}$
    & \cite{HuWang2010} \\[3pt]

    $C_n$, $n\geq 3$
    & $E_{n-2,(n-1)'}e_n
       =(rs)^2e_nE_{n-2,(n-1)'}$
    & \cite{Chen Phd} \\[3pt]

    $D_n$, $n\geq 4$
    & $e_{n-1}e_n=rs\,e_ne_{n-1}$
    & \cite{BGH 2006} \\[3pt]

    $E_6$, $E_7$, $E_8$
    & $\mathcal{E}_{243}e_4=rs\,e_4\mathcal{E}_{243}$
    & \cite{Bai} \\[3pt]

    $F_4$
    & $E_{12332}E_3=(rs)^2E_3E_{12332}$
    & \cite{CHW F4} \\[3pt]

    \hline

    $A_n^{(1)}$, $n\geq 3$
    & $E_{\beta}E_{\eta}=rs\,E_{\eta}E_{\beta}$
    & \cite{HRZ 2008} \\[5pt]

    $\begin{gathered}
      B_n^{(1)},\ C_n^{(1)},\ D_n^{(1)},\\
      E_6^{(1)},\ F_4^{(1)},\ G_2^{(1)}
    \end{gathered}$
    & $e_i e_j=rs\,e_j e_i$
    & \cite{HuZhang2014, Zhang phd 2007} \\[6pt]

    \hline
  \end{tabular}

  \par\medskip
  \begin{minipage}{0.95\textwidth}
    \small
    For finite type $A_n$, set
    $\mathcal{E}_{2,1}=e_2e_1-r^{-1}e_1e_2$
    and $Y=e_2e_3-s e_3e_2$.
    For type $E$, write $\mathcal{E}_{243}=[E_{243}]$
    in the notation of \cite{Bai}.
    For type $F_4$, the simple root generators are denoted by $E_i$.

    \smallskip
    For $A_m^{(1)}$, let
    $\beta=\alpha_m+\alpha_0$ and $\eta=\alpha_1+\alpha_0$, and set
    \[
      E_{\beta}=e_me_0-r e_0e_m,\qquad
      E_{\eta}=e_1e_0-s^{-1}e_0e_1.
    \]

    \smallskip
    In the last row, take $(i,j)=(0,1)$ for
    $B_n^{(1)}$, $D_n^{(1)}$, and $E_6^{(1)}$;
    $(n,0)$ for $C_n^{(1)}$;
    $(0,2)$ for $F_4^{(1)}$; and
    $(2,0)$ for $G_2^{(1)}$.
    Here $n\geq3$ for $B_n^{(1)}$ and $C_n^{(1)}$,
    and $n\geq4$ for $D_n^{(1)}$.
    Further details on the root-vector choices and calculations
    are given in Appendix~\ref{app type D examples}.
  \end{minipage}
  \endgroup
\end{table}

This change in the commutation relation has concrete consequences
for the representation theory of the subalgebras considered here.
We note that a subalgebra generated by two quasi-commutative (resp. commutative) quantum root vectors admits different families of
simple modules: in the quasi-commutative case, both one-dimensional and infinite-dimensional simple modules occur,
whereas in the commutative case only one-dimensional simple modules occur, see Theorem \ref{thm type D isoclass}.
This representation-theoretic difference suggests that the corresponding module categories should be distinct.
To prove that these categories are not equivalent, we use the fact that the center of an algebra is a Morita invariant.
The two subalgebras have non-isomorphic centers, so their categories of all left modules are not equivalent, see Proposition \ref{prop typeD2}.

From the representation-theoretic point of view,
we consider the induction modules from the subalgebra $\mathcal B_{r,s}$ to the whole two-parameter quantum algebra $U_{r,s}(\mathfrak{so}_{2n})$,
and similarly from $\mathcal B_1$ to $U_{q, q^{-1}}(\mathfrak{so}_{2n})$,
where $\mathcal{B}_{r, s}:=\langle e_{n-1}, e_n\rangle \subset U_{r, s}(\mathfrak{so}_{2n})$,
$\mathcal{B}_1:=\langle e_{n-1}, e_n\rangle \subset U_{q, q^{-1}}(\mathfrak{so}_{2n})$.
Frobenius reciprocity then allows us to detect a difference between modules induced from one-dimensional and
infinite-dimensional simple $\mathcal B_{r,s}$-modules.
Namely, if $V_1$ is one-dimensional and $V_2$ is infinite-dimensional, then
$$\operatorname{Hom}_{\mathcal A_{r,s}}
\left(
\mathcal A_{r,s}\otimes_{\mathcal B_{r,s}}V_1, \mathcal A_{r,s}\otimes_{\mathcal B_{r,s}}V_2
\right)=0,$$
see Theorem \ref{thm diff factors}.
The occurrence of infinite-dimensional simple $\mathcal{B}_{r,s}$-modules motivates us to consider
representations beyond finite-dimensional type $1$ weight modules,
whose categories are equivalent in the two-parameter and one-parameter settings.

\textbf{1.3}
This paper is organized as follows.
In section \ref{sec pre}, we fix the root-theoretic notation and recall the Drinfeld--Jimbo presentation of $U_{r,s}(\widehat{\mathfrak g})$,
together with its finite-type subalgebra $U_{r,s}(\mathfrak g)$.
We also introduce the root and weight lattices and the associated group algebras used later.
In section \ref{ssec hopalgdist}, we probe differences in the Hopf algebraic structure between two-parameter (affine) quantum algebra and the one-parameter counterpart.
In section \ref{sec autgps}, we first calculate the automorphism groups of the positive part of two-parameter quantum groups,
via Yakimov's result.
Then we introduce the notion of extended Borel subalgebras and study their automorphisms.
In particular, we obtain several explicit automorphisms in both the root-type and weight-type cases.
In section \ref{sec repdis}, we propose the ``type $D$ phenomenon''.
We prove that $\mathcal B_p\text{-}\mathrm{Mod}$ and $\mathcal B_1\text{-}\mathrm{Mod}$ are not equivalent as abelian categories.
In Section \ref{sec induced}, we investigate the differences between the two-parameter and one-parameter quantum groups from a representation-theoretic perspective.

In Appendix \ref{appA}, we recall the Lyndon bases of type $D_n$,
 including the corresponding orders for the positive and negative parts.
In Appendix \ref{appB}, we calculate triangular relations needed  in Section \ref{sec induced}.

\section{Preliminaries} \label{sec pre}

Let $\mathbb K=\mathbb Q(r,s)\subset\mathbb{C}$, where $r$ and $s$ are algebraically independent elements.
In this paper, all algebras we consider are defined over $\mathbb{C}$.

Let $\mathfrak g$ be a finite-dimensional simple Lie algebra of rank $n$,
with root system $\Phi$ in a Euclidean space $\mathbb R^n$ and with a fixed base
$\Pi=\{\alpha_1,\ldots,\alpha_n\}$ of simple roots.
Let $(\, ,\,)$ denote the inner product on $\mathbb R^n$.
We write $\theta$ for the highest root and $\delta$ for the null root of the corresponding untwisted affine root system.
Set $\alpha_0=\delta-\theta, J=\{1,\ldots,n\}, J_0=J\cup\{0\}.$
For $i\in J_0$, we put $r_i=r^{(\alpha_i,\alpha_i)}, s_i=s^{(\alpha_i,\alpha_i)}.$

When $\mathfrak g$ is of type $D_n$, we use the following realization.
Let $\varepsilon_1,\ldots,\varepsilon_n$ be an orthonormal basis of $\mathbb R^n$.
Then $\Phi(D_n)=
  \{\pm\varepsilon_i\pm\varepsilon_j \mid 1\leq i<j\leq n\}$,
and we choose the base
$\Pi(D_n)=
  \{\alpha_i=\varepsilon_i-\varepsilon_{i+1}\mid 1\leq i<n\}
  \cup
  \{\alpha_n=\varepsilon_{n-1}+\varepsilon_n\}.$

Finally, let
$W=\{\omega_i\mid i\in J_0\}, W'=\{\omega_i'\mid i\in J_0\}$.
We denote by $\Omega=\bigl(\langle \omega_i',\omega_j\rangle\bigr)_{i,j\in J_0}$
the structure constants matrix appearing in the Drinfeld-Jimbo presentation of $U_{r,s}(\widehat{\mathfrak g})$ below.
(The definition of structure constants matrix can be seen in \cite{HuWang2010}).

Write
$\delta=\sum_{i\in J_0}a_i\alpha_i$,
where the $a_i$ are the coefficients of the null root.
For $\mu=\sum_{i\in J_0}m_i\alpha_i$ in the affine root lattice,
we use the notation
\[
\omega_\mu:=\prod_{i\in J_0}\omega_i^{m_i},
\qquad
\omega'_\mu:=\prod_{i\in J_0}(\omega'_i)^{m_i}.
\]
In particular,
\[
\omega_\delta=\prod_{i\in J_0}\omega_i^{a_i},
\qquad
\omega'_\delta=\prod_{i\in J_0}(\omega'_i)^{a_i}.
\]

\begin{defi}\cite{HuZhang2014, Zhang phd 2007}\label{def DJ}
  The unital associative algebra $U_{r,s}(\widehat{\mathfrak{g}})$ over $\mathbb{C}$ is generated by the elements $e_j, f_j, \omega_j^{\pm}, (\omega'_j)^{\pm}, \gamma^{\pm \frac{1}{2}}, (\gamma')^{\pm \frac{1}{2}}, D^{\pm}, (D')^{\pm}$, where $j\in J_0$, and satisfies the following relations:

  \medskip
  $(\textrm{1})$
  $\gamma^{\pm\frac{1}{2}}, (\gamma')^{\pm\frac{1}{2}}$ are central with $\gamma = (\omega'_\delta)^{-1}$, $\gamma' = (\omega_\delta)^{-1}$, $\gamma\gamma' = (rs)^c$ (where $c$ is the canonical central element of $\widehat{\mathfrak{g}}$),
   such that $\omega_i \omega_i^{-1} = \omega'_i(\omega'_i)^{-1} = 1 = DD^{-1} = D'D'^{-1}$, and
  \begin{align*}
    [\omega_i^{\pm 1}, \omega_j^{\pm 1}] &= [\omega_i^{\pm1}, D^{\pm1}] = [(\omega'_j)^{\pm1}, D^{\pm1}] = [\omega_i^{\pm1}, (D')^{\pm1}] = 0 \\
    &= [\omega_i^{\pm 1}, (\omega'_j)^{\pm 1}] = [(\omega'_j)^{\pm1}, (D')^{\pm1}] = [(D')^{\pm1}, D^{\pm1}] = [(\omega'_i)^{\pm 1}, (\omega'_j)^{\pm 1}].
  \end{align*}

  $(\textrm{2})$ For $i, j\in J_0$,
  \begin{align*}
    D e_i D^{-1} &= r_i^{\delta_{0i}} e_i, \qquad\qquad D f_i D^{-1} = r_i^{-\delta_{0i}} f_i, \\
    \omega_j e_i (\omega_j)^{-1} &= \langle \omega'_i, \omega_j \rangle e_i, \qquad\quad \omega_j f_i (\omega_j)^{-1} = \langle \omega'_i, \omega_j \rangle^{-1} f_i.
  \end{align*}

  $(\textrm{3})$ For $i, j\in J_0$,
  \begin{align*}
    D' e_i (D')^{-1} &= s_i^{\delta_{0i}} e_i, \qquad\qquad D' f_i (D')^{-1} = s_i^{-\delta_{0i}} f_i, \\
    \omega'_j e_i (\omega'_j)^{-1} &= \langle \omega'_j, \omega_i \rangle^{-1} e_i, \qquad\quad \omega'_j f_i (\omega'_j)^{-1} = \langle \omega'_j, \omega_i \rangle f_i.
  \end{align*}

  $(\textrm{4})$ For $i, j\in J_0$, we have
  $$[\,e_i, f_j\,] = \frac{\delta_{ij}}{r_i - s_i}(\omega_i - \omega'_i).$$

  $(\textrm{5})$ For any $i\ne j$, we have the $(r,s)$-Serre relations
  \begin{gather*}
    \bigl(\operatorname{ad}_l e_i\bigr)^{1-a_{ij}}(e_j) = 0, \\
    \bigl(\operatorname{ad}_r f_i\bigr)^{1-a_{ij}}(f_j) = 0,
  \end{gather*}
  where the $(a_{ij})$ is the Cartan matrix of $\widehat{\mathfrak{g}}$,
the left-adjoint action $\operatorname{ad}_l e_i$ and the right-adjoint action $\operatorname{ad}_r f_i$ are defined as follows:
  $$
  \operatorname{ad}_l a (b) = \sum_{(a)} a_{(1)} b S(a_{(2)}), \quad
  \operatorname{ad}_r a (b) = \sum_{(a)} S(a_{(1)}) b a_{(2)}, \quad
  \forall\; a, b\in U_{r,s}(\widehat{\mathfrak{g}}),
  $$
  and the coproduct $\Delta(a) = \sum_{(a)} a_{(1)} \otimes a_{(2)}$ is given by
  \begin{gather*}
    \Delta(\omega_i^{\pm1}) = \omega_i^{\pm1} \otimes \omega_i^{\pm1}, \qquad
    \Delta((\omega'_i)^{\pm1}) = (\omega'_i)^{\pm1} \otimes (\omega'_i)^{\pm1}, \\
    \Delta(e_i) = e_i \otimes 1 + \omega_i \otimes e_i, \qquad \Delta(f_i) = 1 \otimes f_i + f_i \otimes \omega'_i, \\
    \varepsilon(\omega_i^{\pm}) = \varepsilon((\omega'_i)^{\pm1}) = 1, \qquad
    \varepsilon(e_i) = \varepsilon(f_i) = 0, \\
    S(\omega_i^{\pm1}) = \omega_i^{\mp1}, \qquad
    S((\omega'_i)^{\pm1}) = (\omega'_i)^{\mp1}, \\
    S(e_i) = -\omega_i^{-1}e_i, \qquad S(f_i) = -f_i (\omega'_i)^{-1}.
  \end{gather*}
\end{defi}

\begin{remark}\label{rmk finaff}
  If we restrict the generators on $\{e_j, f_j, \omega_j, \omega_j'\mid j\in J\}$, we will obtain the definition of $U_{r,s}(\mathfrak{g})$, see \cite{BGH 2006}.
\end{remark}

For later use, we also fix the following lattice notation for the underlying
finite root system. Let \(Q=\Lambda_R=\bigoplus_{i\in I}\mathbb Z\alpha_i\)
and \(P=\Lambda_W=\bigoplus_{i\in I}\mathbb Z\lambda_i\) be respectively the
root lattice and the weight lattice of \(\mathfrak g\), where \(\lambda_i\) are
the fundamental weights. More generally, if \(L\) is a lattice with
\(Q\subseteq L\subseteq P\), we write \(\mathbb K[L\oplus L]\) for the group
algebra generated by \(\omega_\lambda^{\pm1}\) and
\((\omega'_\lambda)^{\pm1}\), \(\lambda\in L\), with
\(\omega_\lambda\omega_\mu=\omega_{\lambda+\mu}\),
\(\omega'_\lambda\omega'_\mu=\omega'_{\lambda+\mu}\), and
\(\omega_\lambda\omega'_\mu=\omega'_\mu\omega_\lambda\).

\section{Two- vs. one-parameter Hopf algebraic distinctions}\label{ssec hopalgdist}

In this section, we will probe differences in the Hopf algebraic structure between two-parameter (affine) quantum algebra and the one-parameter counterpart.
Assume that $C$ is any coalgebra and $g, h$ are group-like elements of $C$,
then the $(g, h)$-skew primitive element $x\in C$ is defined by $\Delta(x)=x\otimes g+ h\otimes x$.
The set of all $(g, h)$-skew primitive elements in $C$ is denoted by $P_{g, h}(C)$.
Firstly we need the following lemma.

\begin{lemm}\label{lemm skewpri}
We have the following results concerning the skew-primitive elements of $U_{r, s}(\mathfrak{g})$:

$(1)$\quad
$P_{1,\omega_i}=\mathbb{C}(1-\omega_i)+\mathbb{C}e_i$,
\quad $
P_{1,\omega_i'^{-1}}=\mathbb{C}(1-\omega_i'^{-1})+\mathbb{C}f_i\omega_i'^{-1}$,

$\qquad \ P_{1,\sigma}=\mathbb{C}(1-\sigma)
\quad \textrm{if}\quad
 \sigma\not\in \{\omega_j, \omega_j'^{-1}\mid j\in J\};$

$(2)\quad
P_{\omega'_i,1}=\mathbb{C}(1-\omega'_i)+\mathbb{C}f_i$,
\quad $
P_{\omega_i^{-1},1}=\mathbb{C}(1-\omega_i^{-1})+\mathbb{C}e_i\omega_i^{-1}$,

$\qquad \ P_{\sigma,1}=\mathbb{C}(1-\sigma)
\quad \textrm{if}\quad
 \sigma\not\in\{ \omega'_j, \omega_j^{-1} \mid j\in J\}.$
\end{lemm}

\begin{proof}
Similar to that of Lemma 4.3 in \cite{HuWang2009}.
\end{proof}

\begin{theorem}\label{thm isofin}
Let $\mathfrak g$ be a finite-dimensional complex simple
Lie algebra, and let $q\in\mathbb C^\times$ with
$q\ne\pm1$.
If $\mathfrak g\not\cong\mathfrak{sl}_2$, then there exists
a $\mathbb C$-Hopf algebra isomorphism
\[
\phi:U_{r,s}(\mathfrak g)
\longrightarrow U_{q,q^{-1}}(\mathfrak g)
\]
if and only if
$(r,s)=(q,q^{-1})\,\text{or} \,(r,s)=(q^{-1},q)$.

If $\mathfrak g\cong\mathfrak{sl}_2$, then
$U_{r,s}(\mathfrak{sl}_2)$ and
$U_{q,q^{-1}}(\mathfrak{sl}_2)$ are isomorphic
as $\mathbb C$-Hopf algebras if and only if
$q^2=rs^{-1}\;\text{or}\;q^2=r^{-1}s.$
\end{theorem}

\begin{proof}
  Write $\tilde{e}_i, \tilde{f}_i, \tilde{\omega}_i, \tilde{\omega}'_i$ to distinguish the generators of $U_{q, q^{-1}}(\mathfrak{g})$.
  For any $i\in J$, we have
  \begin{align*}
    \Delta(\phi(e_i))\; &=(\phi\otimes\phi)(\Delta(e_i))=(\phi\otimes\phi)(e_i\otimes 1+\omega_i\otimes e_i) \\
     & =\phi(e_i)\otimes 1+\phi(\omega_i)\otimes\phi(e_i).
  \end{align*}
  Notice that $\phi(\omega_i)$ is a group-like element,
  which means $\phi(e_i)$ is a $\bigl(1, \phi(\omega_i)\bigr)$-skew primitive element.
  Since $\phi$ is surjective, then $\phi(\omega_i)=\tilde{\omega}_{j_i}$ or $\tilde{\omega}_{j_i}'^{-1}$, where $\tilde{\omega}_{j_i}$, $\tilde{\omega}_{j_i}'^{-1}$ are some group-like elements of $U_{q, q^{-1}}(\mathfrak{g})$.
  By Lemma \ref{lemm skewpri},
  If $\phi(\omega_i)=\tilde{\omega}_{j_i}$, we may assume
  $$\phi(e_i)=a_i\tilde{e}_{j_i}+b_i (1-\tilde{\omega}_{j_i}), \quad a_i, b_i\in \mathbb{C}.$$
  Substituting these into $\phi(\omega_ie_i)=\langle \omega'_i, \omega_i \rangle \phi( e_i\omega_i), 1\leq i\leq n,$
  we obtain
  \begin{equation}\label{equ Kji}
    \tilde{\omega}_{j_i}\bigl(a_i \tilde{e}_{j_i}+b_i(1-\tilde{\omega}_{j_i})\bigr)=\langle \omega'_i, \omega_i \rangle
    \bigl(a_i \tilde{e}_{j_i}+b_i(1-\tilde{\omega}_{j_i})\bigr) \tilde{\omega}_{j_i}.
  \end{equation}
  Since $\langle \omega'_i, \omega_i \rangle\neq 1$, it follows immediately that $b_i=0$.

  If $\phi(\omega_i)=\tilde{\omega}_{j_i}'^{-1}$,
   we may assume
  $$\phi(e_i)=c_i \tilde{f}_{j_i}\tilde{\omega}_{j_i}'^{-1}+d_i(1-\tilde{\omega}_{j_i}'^{-1}), \quad c_i, d_i\in \mathbb{C}.$$
  Analogously, substituting into the same relation yields
  \begin{equation}\label{equ Kji2}
    \tilde{\omega}_{j_i}'^{-1}\bigl(c_i \tilde{f}_{j_i}\tilde{\omega}_{j_i}'^{-1}+d_i(1-\tilde{\omega}_{j_i}'^{-1})\bigr)=\langle \omega'_i, \omega_i \rangle\bigl(c_i \tilde{f}_{j_i}\tilde{\omega}_{j_i}'^{-1}+d_i(1-\tilde{\omega}_{j_i}'^{-1})\bigr) \tilde{\omega}_{j_i}'^{-1},
  \end{equation}
  and thus $d_i=0$.
  According to the diagram automorphism group of $\frak{g}$,
  we need to consider the following three situations:

  (I) $\frak{g}=B_n, C_n, E_7, E_8, F_4, G_2$. Applying (\ref{equ Kji}) or (\ref{equ Kji2}), we can check  that $\phi$ must preserve the length of root vectors.
  Notice that $\phi([e_i, e_j])=[\phi(e_i), \phi(e_j)]$,
  we conclude that $\phi$ must preserve the adjacency relations of Dynkin diagram.
  Specifically, if the simple roots corresponding to \(e_i\) and \(e_j\) are connected in the Dynkin diagram,
  then the simple roots corresponding to \(\phi(e_i)\) and \(\phi(e_j)\) must also be connected.
  Therefore $j_i=i\in J$ in these cases.
  Since $\phi$ preserves Serre relations, we deduce that either $\phi(e_i)=a_i \tilde{e}_i$ for all $i$, or $\phi(e_i)=b_i \tilde{f}_{i}\tilde{\omega}_{i}'^{-1}$ for all $i$.
   If $\phi(e_i)=a_i \tilde{e}_i$, then substituting into
   \begin{equation}\label{}
     \phi(\omega_ie_j)=\langle \omega'_j, \omega_i \rangle \phi( e_j\omega_i),\quad i, j\in J,
   \end{equation}
   yields $\langle \tilde{\omega}_{j}', \tilde{\omega}_{i}\rangle=\langle \omega'_j, \omega_i \rangle$.
   This implies $(r, s)=(q, q^{-1})$.
   If $\phi(e_i)=b_i \tilde{f}_{i}\tilde{\omega}_{i}'^{-1}$, then $\langle \tilde{\omega}_{i}', \tilde{\omega}_{j}\rangle^{-1}=\langle \omega'_j, \omega_i \rangle$, which implies  $(r, s)=(q^{-1}, q)$.

  (II) $\frak{g}= A_n\, (n>1), \, D_n\, (n\geq 4)$ or $E_6$.
  Here it suffices to prove the case of $\frak{g}= A_2$.
  Consider the following equations:
  \begin{align}
     & \phi(\omega_1e_1)=\langle \omega'_1, \omega_1 \rangle \phi( e_1\omega_1), \label{equ 1}\\
     & \phi(\omega_1e_2)=\langle \omega'_2, \omega_1 \rangle \phi( e_2\omega_1), \label{equ 2}
  \end{align}
  Similarly, we need to consider the following two situations:

  If $\phi$ preserves positive root vectors:
  From (\ref{equ 1}), we obtain $rs^{-1}=q^2$.
  From (\ref{equ 2}), it follows that $q^{-1}=s$.

  If $\phi$ maps positive root vectors to negative ones:
  From (\ref{equ 1}), we have $rs^{-1}=q^{-2}$.
  From (\ref{equ 2}), it follows that $q=s$.

  In either situation,  we conclude that either $(r, s)=(q, q^{-1})$, or $(r, s)=(q^{-1}, q)$.

 (III) The $\mathfrak{sl}_2$ case can be found in \cite{BenkartW2004}.
\end{proof}

\begin{remark}
   Benkart-Witherspoon gave a proof for type $A$ in \cite{BenkartW2004},
where the symmetric case: $(r, s)=(q^{-1}, q)$ was not considered.
\end{remark}

One can extend this theorem to the affine case easily, since both $U_{r, s}(\widehat{\mathfrak{g}})$ and $U_{q, q^{-1}}(\widehat{\mathfrak{g}})$ are pointed Hopf algebras.

\begin{coro}\label{cor iso}
   Assume that $\mathfrak{g}\neq\mathfrak{sl}_2$.
   There exists a $\mathbb C$-Hopf algebra isomorphism
   $$\phi: U_{r, s}(\widehat{\frak{g}})\longrightarrow U_{q, q^{-1}}(\widehat{\frak{g}})$$
   for  $q\ne \pm1$, if and only if $(r, s)=(q, q^{-1})$, or $(r, s)=(q^{-1}, q)$.
\end{coro}

\begin{proof}
Let $\phi:U_{r,s}(\widehat{\mathfrak g})\longrightarrow U_{q,q^{-1}}(\widehat{\mathfrak g})$ be a $\mathbb C$-Hopf algebra isomorphism.
Write $\tilde{e}_i, \tilde{f}_i, \tilde{\omega}_i, \tilde{\omega}'_i\, (i\in J_0)$ to distinguish the generators of $U_{q, q^{-1}}(\widehat{\mathfrak{g}})$.
Denote by $U_{r,s}(\mathfrak g)\subset U_{r,s}(\widehat{\mathfrak g})$ the Hopf subalgebra generated by $e_i,f_i,\omega_i^{\pm1},(\omega_i')^{\pm1}$, $i\in J$.
This is the finite-type two-parameter quantum group associated with $\mathfrak g$.

Since both Hopf algebras are pointed, $\phi$ sends group-like elements to group-like elements.
Moreover, $\Delta(\phi(e_i))=(\phi\otimes\phi)\Delta(e_i)=\phi(e_i)\otimes 1+\phi(\omega_i)\otimes \phi(e_i)$,
so $\phi(e_i)$ is a $(1,\phi(\omega_i))$-skew primitive element.
By the same description of skew primitive elements used in Lemma \ref{lemm skewpri} and Theorem \ref{thm isofin},
there is a permutation $\sigma:J_0\to J_0$ such that, for each $i\in J_0$,
the image of the rank-one part generated by $e_i,f_i,\omega_i^{\pm1},(\omega_i')^{\pm1}$ is contained in the
rank-one part of $U_{q,q^{-1}}(\widehat{\mathfrak g})$ indexed by $\sigma(i)$.
More precisely, up to nonzero scalar multiples, $\phi(e_i)$ is either $\widetilde e_{\sigma(i)}$ or $\widetilde f_{\sigma(i)}\widetilde\omega_{\sigma(i)}^{\prime\,-1}$,
and similarly for $\phi(f_i)$.

The quantum Serre relations force $\sigma$ to preserve the adjacency relations of the affine Dynkin diagram.
Hence $\sigma$ is an automorphism of the affine Dynkin diagram.
In particular, the subdiagram with vertex set $\sigma(J)$ is isomorphic to the finite Dynkin diagram of $\mathfrak g$.
Therefore the Hopf subalgebra of $U_{q,q^{-1}}(\widehat{\mathfrak g})$ generated by
$\widetilde e_j,\widetilde f_j,\widetilde\omega_j^{\pm1},
(\widetilde\omega_j')^{\pm1}$, $j\in \sigma(J)$,
is isomorphic, after relabelling the simple roots, to $U_{q,q^{-1}}(\mathfrak g)$.

Thus the restriction of $\phi$ to $U_{r,s}(\mathfrak g)$ gives a finite-type Hopf algebra isomorphism
\[
  \phi|_{U_{r,s}(\mathfrak g)}:
  U_{r,s}(\mathfrak g)\longrightarrow U_{q,q^{-1}}(\mathfrak g),
\]
up to a Dynkin diagram relabelling in the target.
Applying Theorem \ref{thm isofin} and using the assumption $\mathfrak g\neq\mathfrak{sl}_2$,
we obtain $(r,s)=(q,q^{-1})$ or $(r,s)=(q^{-1},q)$.
This proves this Corollary.
\end{proof}

\section{The automorphism of the extended Borel subalgebra}\label{sec autgps}
In this section, we turn to the study of algebra automorphisms of two-parameter quantum groups.
Using Yakimov's rigidity theorem, we first determine the automorphism groups of their positive parts.
We then introduce the notion of extended Borel subalgebras (see Definition \ref{def extended}) and describe several explicit automorphisms of these algebras, in types $A_n$ and $D_4$.

In the Dynkin diagram of type $D_4$, the simple roots $\alpha_1,\alpha_3,\alpha_4$ play symmetric roles
and the diagram automorphism group is isomorphic to $S_3$.
However, for the root-type extended Borel algebra, only the transposition exchanging $\alpha_3$ and $\alpha_4$ can be lifted to an algebra automorphism.
The Dynkin diagram automorphisms interchanging $\alpha_1$ with $\alpha_3$ or $\alpha_4$ cannot be realized at the level of the root lattice.
To recover the full symmetry, we enlarge the root lattice to the weight lattice and consider the corresponding
weight-type extended Borel algebra.

\begin{theorem}{(Rigidity Theorem)}\label{Thm Yakimov}\cite{Yakimov 2014}
  For all simple Lie algebras $\mathfrak{g}$ of rank $n>1$,
base fields $\mathbb{C}$, and deformation parameters $q\in \mathbb{C}$  which are not roots of unity,
$$Aut(U_q^+(\mathfrak{g}))\cong \mathbb{T}^n \rtimes Aut(\Gamma),$$
where $U_q^+(\mathfrak{g})$ is the subalgebra generated by $e_i\, (i\in J)$,
$\mathbb{T}^n=(\mathbb{C}^\ast)^n$ is the torus, Aut$(\Gamma)$ is the automorphism group of the Dynkin diagram of $\mathfrak{g}$.
\end{theorem}

Assume that $Q$ is the root lattice of $\mathfrak{g}$,
and that $U_q^+(\mathfrak{g})$ is a $\mathbb{C}$-algebra graded by $Q$.
Given a $2$-cocycle $\mathbf{p}\in Z^2(Q,\mathbb{C}^\ast)$,
we define a new algebra $U_{q,\mathbf{p}}^+(\mathfrak{g})$
by twisting the multiplication:
$$
E_{\alpha_1}\ast E_{\alpha_2}
=\mathbf{p}(\alpha_1,\alpha_2)E_{\alpha_1}E_{\alpha_2},
\quad
\alpha_i\in Q,\;
E_{\alpha_i}\in U_q^+(\mathfrak{g})_{\alpha_i},
\quad i=1,2.
$$
Define
$$
\mathbf{r}:Q\times Q\rightarrow\mathbb{C}^\ast,
\quad
\mathbf{r}(\alpha_1,\alpha_2)
=\mathbf{p}(\alpha_1,\alpha_2)
 \mathbf{p}(\alpha_2,\alpha_1)^{-1},
\quad \alpha_1,\alpha_2\in Q.
$$
Denote
$$
\operatorname{Aut}(\Gamma,\mathbf{p})
=
\left\{
\theta\in\operatorname{Aut}(\Gamma)
\mid
\mathbf{r}(\theta(\alpha),\theta(\alpha'))
=\mathbf{r}(\alpha,\alpha'),
\quad \forall\alpha,\alpha'\in\Pi
\right\}.
$$
In Propositions \ref{aut autA2+}, \ref{aut D4}, and \ref{aut E6},
we shall see that this subgroup is trivial in the
two-parameter cases considered there.

\begin{theorem} \cite{Yakimov 2014}
For every simple Lie algebra $\mathfrak{g}$ of rank $\ell>1$,
$q\in\mathbb{C}^\ast$, and a torsion-free $2$-cocycle
$\mathbf{p}\in Z^2(Q,\mathbb{C}^\ast)$ satisfying
\begin{equation}\label{Yakimov condition}
q_\alpha\mathbf{r}(\alpha,\alpha')\neq1,
\quad
q_\alpha^{-1}\mathbf{r}(\alpha,\alpha')\neq1,
\quad
\forall\alpha,\alpha'\in\Pi
\text{ such that }a_{\alpha,\alpha'}=-1,
\end{equation}
the map
$$
\Upsilon^+:
\mathbb{T}^{\ell}\rtimes\operatorname{Aut}(\Gamma,\mathbf{p})
\rightarrow
\operatorname{Aut}(U_{q,\mathbf{p}}^+(\mathfrak{g})),
$$
defined on the generators by
$$
\Upsilon^+(t,\theta)(e_i)
=t_{\theta(i)}e_{\theta(i)},
\quad
t=(t_i)_{i\in J}\in\mathbb{T}^{\ell},
\quad
\theta\in\operatorname{Aut}(\Gamma,\mathbf{p}),
\quad i\in J,
$$
is a group isomorphism.
\end{theorem}

In \cite{HuPei2008}, Hu-Pei proved that
\begin{theorem}\label{Thm HuPei}
$U_{r,s}^+(\mathfrak{g})$ and
$U_{r',s'}^+(\mathfrak{g})$ are $2$-cocycle deformations
of each other if $rs^{-1}=r'(s')^{-1}$
or $rs^{-1}=(r')^{-1}s'$.
\end{theorem}

Choose square roots of $r$ and $s$, and assume that
$r^{\frac12}s^{-\frac12}=q$.
We identify $U_q^+(\mathfrak{g})$ with
$U_{q,q^{-1}}^+(\mathfrak{g})$.
In the preceding construction, we now take
$$
\mathbf{p}(\mu,\nu)
=
\left(r^{-\frac12}s^{-\frac12}\right)^{
\langle\mu,\nu\rangle},
\quad \mu,\nu\in Q,
$$
where $\langle\,\cdot,\cdot\,\rangle$ is the bilinear
form defined below.
The symbol $\ast$ denotes the same twisted multiplication
as above, with $xy$ denoting the original multiplication
in $U_{q,q^{-1}}^+(\mathfrak{g})$.
With this choice of $\mathbf{p}$, the assignment
$e_i\mapsto e_i$ gives a graded algebra isomorphism
$$
U_{r,s}^+(\mathfrak{g})
\cong U_{q,\mathbf{p}}^+(\mathfrak{g}),
$$
where the algebra on the right is the cocycle twist of
$U_{q,q^{-1}}^+(\mathfrak{g})$; see \cite{HuPei2008}.

The bilinear form $\langle\,\cdot,\cdot\,\rangle$
is defined on the root lattice $Q$ by
$$
\langle i,j\rangle
:=\langle\alpha_i,\alpha_j\rangle
=
\left\{
\begin{array}{ll}
d_i a_{ij}, & \hbox{$i<j$;} \\
d_i,       & \hbox{$i=j$;} \\
0,         & \hbox{$i>j$,}
\end{array}
\right.
$$
where $(a_{ij})$ is the Cartan matrix of finite type,
and $\{d_i\mid i\in J\}$ is the set of relatively prime
positive integers such that
$d_i a_{ij}=d_j a_{ji}$ for $i,j\in J$.

In type $D_n$, we modify the two entries
$\langle n-1,n\rangle$ and $\langle n,n-1\rangle$:
\begin{equation}\label{equ redefi}
\langle n-1,n\rangle=-1,
\quad
\langle n,n-1\rangle=1.
\end{equation}

According to \cite{BGH 2006},  two-parameter quantum groups do not admit the usual Lusztig symmetries when the rank $n>2$.
This indicates that the structure of $U_{r,s}(\mathfrak g)$ is considerably more rigid than that of the one-parameter quantum group $U_q(\mathfrak g)$.

\iffalse
\begin{prop}
The positive part \(U_{r,s}^{+}(\mathfrak{sl}_n)\) is a domain.
\end{prop}

\begin{proof}
By Proposition 3.4 of Benkart--Kang--Lee, the algebra
\(U_{r,s}^{+}(\mathfrak{sl}_n)\) is an iterated skew polynomial algebra over
\(\mathbb K\). More precisely, with respect to the root vectors
\(\mathcal E_{i,j}\), one has
\[
U_{r,s}^{+}(\mathfrak{sl}_n)
 =
\mathbb K[\mathcal E_{1,1}]
[\mathcal E_{2,1};\varphi_{2,1},\vartheta_{2,1}]
\cdots
[\mathcal E_{n-1,n-1};\varphi_{n-1,n-1},\vartheta_{n-1,n-1}],
\]
where each \(\varphi_{i,j}\) is an algebra automorphism and
\(\vartheta_{i,j}\) is a \(\varphi_{i,j}\)-derivation of the subalgebra
generated by the preceding root vectors.

We use the standard fact that if \(R\) is a domain, \(\sigma\) is an injective
endomorphism of \(R\), and \(\delta\) is a \(\sigma\)-derivation, then the Ore
extension \(R[x;\sigma,\delta]\) is again a domain. Indeed, for nonzero
polynomials \(f=a_mx^m+\cdots\) and \(g=b_nx^n+\cdots\), their product has
leading coefficient \(a_m\sigma^m(b_n)\), which is nonzero since \(R\) is a
domain and \(\sigma\) is injective. Hence \(fg\neq 0\).

Starting from the field \(\mathbb K\), which is a domain, and applying this
fact successively to the above iterated skew polynomial presentation, we obtain
that \(U_{r,s}^{+}(\mathfrak{sl}_n)\) is a domain.
\end{proof}
\fi

The next proposition will demonstrate the specifics of the rigidity of $U_{r, s}^+(\mathfrak{g})$.

\begin{prop}\label{aut autA2+}
  For $\mathfrak{g}=\mathfrak{sl}_3$,
  $$\text{Aut}(U_{r, s}^+(\mathfrak{g}))\cong \mathbb{T}^2.$$
  That is to say,  the nontrivial Dynkin diagram automorphism of $\mathfrak{sl}_3$ does not induce an algebraic automorphism of $U_{r, s}^+(\mathfrak{g})$.
\end{prop}

\begin{proof}
  We use Yakimov's result to determine the automorphism group.
    Assume $\sigma$ is the diagram automorphism exchanging $\alpha_1$ with $\alpha_2$.
     By direct calculations, we find that (\ref{Yakimov condition}) is always satisfied, and
$$\mathbf{r}(\alpha_1, \alpha_2)=\mathbf{p}(\alpha_1, \alpha_2)\mathbf{p}(\alpha_2, \alpha_1)^{-1}=r^{\frac{1}{2}}s^{\frac{1}{2}}, \quad
\mathbf{r}(\alpha_2, \alpha_1)=\mathbf{p}(\alpha_2, \alpha_1)\mathbf{p}(\alpha_1, \alpha_2)^{-1}=r^{-\frac{1}{2}}s^{-\frac{1}{2}}.$$
Obviously, $\mathbf{r}(\alpha_1, \alpha_2)\neq \mathbf{r}(\alpha_2, \alpha_1)$ unless $rs=1$.
Then we conclude that $\sigma$ does not induce an automorphism of $\text{Aut}(U_{r, s}^+(\mathfrak{g}))$.
\end{proof}

\begin{remark}
  We also verify directly that the nontrivial Dynkin
diagram automorphism cannot be realized by interchanging
the generators.
Suppose that an algebra automorphism $\phi$ satisfies
$\phi(e_1)=a e_2,\, \phi(e_2)=b e_1, \,a,b\in\mathbb C^\times.$
Under the above assumption, we have
$\phi(\mathcal E_{12})=ab(e_2e_1-se_1e_2)$.
Since $e_1\mathcal E_{12}=r\mathcal E_{12}e_1$,
applying $\phi$ and dividing by $a^2b$ gives
$e_2(e_2e_1-se_1e_2)=r(e_2e_1-se_1e_2)e_2.$
Expanding this identity yields
$e_2^2e_1-(r+s)e_2e_1e_2+rs\,e_1e_2^2=0.$
  Equivalently, we have
  \begin{equation}\label{equ serre1}
    e_1e_2^2-(r^{-1}+s^{-1})e_2e_1e_2+r^{-1}s^{-1}e_2^2e_1=0.
  \end{equation}
From the Serre relations of $U_{r, s}(\mathfrak{sl}_3)$, we have
\begin{equation}\label{equ serre2}
  e_1e_2^2-(r+s)e_2e_1e_2+rse_2^2e_1=0.
\end{equation}
Subtracting (\ref{equ serre1}) from (\ref{equ serre2}), we have
$$e_2\Bigl[(rs-r^{-1}s^{-1})e_2e_1-(r+s-r^{-1}-s^{-1})e_1e_2\Bigr]=0.$$
By Proposition 3.4 in \cite{Benkart06}, $U_{r,s}^+(\mathfrak{sl}_3)$ admits an iterated skew polynomial structure.
It follows that this algebra is a domain and hence contains no zero divisors.
Consequently, we obtain
$rs - r^{-1}s^{-1} = r + s - r^{-1} - s^{-1} = 0.$
This forces $rs=1$, yielding a contradiction.
\end{remark}

\begin{prop}\label{aut D4}
   For $\mathfrak{g}=\mathfrak{so}_8$,
  $$\text{Aut}(U_{r, s}^+(\mathfrak{g}))\cong \mathbb{T}^4.$$
\end{prop}

\begin{proof}
  We still use Yakimov's result  to determine the automorphism group.
  We first remark that Theorem \ref{Thm HuPei} is still true for type $D_4$, via the redefinition (\ref{equ redefi}).
Indeed,
\begin{align*}
  e_3\ast e_4-rs e_4\ast e_3 & =\,\mathbf{p}(\alpha_3, \alpha_4)e_3e_4-rs\mathbf{p}(\alpha_4, \alpha_3)e_4e_3 \\
   & =\,r^{\frac{1}{2}}s^{\frac{1}{2}}e_3e_4-rsr^{-\frac{1}{2}}s^{-\frac{1}{2}}e_4e_3  \\
   & =\,r^{\frac{1}{2}}s^{\frac{1}{2}}(e_3e_4-e_4e_3)\\
   &=\, 0.
\end{align*}
Assume $\sigma$ is the diagram automorphism of $\mathfrak{so}_8$ exchanging $\alpha_3$ with $\alpha_4$.
We also know from the verifications that
$$rs=\mathbf{r}(\alpha_3, \alpha_4)\neq \mathbf{r}(\alpha_4, \alpha_3)=r^{-1}s^{-1}, \quad \text{unless}\; rs=1.$$
This means that $\sigma$ does not induce an automorphism of $\text{Aut}(U_{r, s}(\mathfrak{n}_+))$.

Similarly, assume that $\sigma'$ is the diagram automorphism of $\mathfrak{so}_8$ exchanging $\alpha_1$ with $\alpha_3$.
Then we calculate that
$$\mathbf{r}(\alpha_1, \alpha_2)=\mathbf{p}(\alpha_1, \alpha_2)\mathbf{p}(\alpha_2, \alpha_1)^{-1}=r^{\frac{1}{2}}s^{\frac{1}{2}}, \quad
\mathbf{r}(\alpha_3, \alpha_2)=\mathbf{p}(\alpha_3, \alpha_2)\mathbf{p}(\alpha_2, \alpha_3)^{-1}=r^{-\frac{1}{2}}s^{-\frac{1}{2}}.$$
Then $\mathbf{r}(\alpha_1, \alpha_2)\neq \mathbf{r}(\alpha_3, \alpha_2)$ unless $rs=1$,
which means that $\sigma'$ does not induce an automorphism of $\text{Aut}(U_{r, s}(\mathfrak{n}_+))$.

By the same token, we can verify that any nontrivial automorphism of $\mathfrak{so}_8$ does not induce an automorphism of $\text{Aut}(U_{r, s}^+(\mathfrak{g}))$.
\end{proof}

\begin{prop}\label{aut E6}
     For $\mathfrak{g}=E_6$,
  $$\text{Aut}(U_{r, s}^+(\mathfrak{g}))\cong \mathbb{T}^6.$$
\end{prop}

\begin{proof}
    We use Yakimov's result to determine the automorphism group.
    Assume $\sigma$ is the diagram automorphism exchanging $\alpha_i$ with $\alpha_{7-i}, \, i=1, 3$.
By direct calculations, we find that (\ref{Yakimov condition}) is always satisfied, and
$$\mathbf{r}(\alpha_3, \alpha_4)=\mathbf{p}(\alpha_3, \alpha_4)\mathbf{p}(\alpha_4, \alpha_3)^{-1}=r^{\frac{1}{2}}s^{\frac{1}{2}}, \quad
\mathbf{r}(\alpha_5, \alpha_4)=\mathbf{p}(\alpha_5, \alpha_4)\mathbf{p}(\alpha_4, \alpha_5)^{-1}=r^{-\frac{1}{2}}s^{-\frac{1}{2}}.$$
Obviously, $\mathbf{r}(\alpha_3, \alpha_4)\neq \mathbf{r}(\alpha_5, \alpha_4)$ unless $rs=1$.
Then we conclude that $\sigma$ does not induce an automorphism of $\text{Aut}(U_{r, s}^+(\mathfrak{g}))$.
\end{proof}

Based on the differences in their automorphism groups, we conclude that
\begin{prop}
For $\mathfrak{g}=\mathfrak{sl}_{n+1}\,(n\geq 2), \mathfrak{so}_{2n}\, (n\geq 4)$ or $E_6$, we have
$$U_{r, s}^+(\mathfrak{g})\ncong U_{q, q^{-1}}^+(\mathfrak{g})$$
as $\mathbb C$-algebras.
\end{prop}

\begin{remark}
  To remedy the lack of symmetry in two-parameter quantum groups, we need to adjoin all group-like elements.
  This motivates our study of extended Borel algebras.
\end{remark}

\begin{defi}\label{def extended}
Let $L$ be either the root lattice $Q$ or the weight lattice $P$.
Use the pairing convention of Definition~2.1 and, when $L=P$, fix a bicharacter extension of that pairing to the two copies of $L$.
Write $\omega_\lambda$ and $\omega'_\lambda$ for the canonical group elements corresponding to $(\lambda,0)$ and $(0,\lambda)$ in $\mathbb K[L\oplus L]$, respectively.
Thus
\[
\omega_\lambda\omega_\mu=\omega_{\lambda+\mu},\qquad
\omega'_\lambda\omega'_\mu=\omega'_{\lambda+\mu},\qquad
\omega_\lambda\omega'_\mu=\omega'_\mu\omega_\lambda,
\]
with $\omega_0=\omega'_0=1$.

The \emph{$L$-type extended Borel algebra} $\mathcal B^L_{r,s}(\mathfrak g)$ is generated by $U^+_{r,s}(\mathfrak g)$ and $\mathbb K[L\oplus L]$, subject to the cross relations
\begin{align*}
&\omega_\lambda e_i\omega_\lambda^{-1}=\langle\omega'_i,\omega_\lambda\rangle e_i, \quad
\omega_\lambda f_i\omega_\lambda^{-1}=\langle\omega'_i,\omega_\lambda\rangle^{-1}f_i,\\
&\omega'_\lambda e_i(\omega'_\lambda)^{-1}
=\langle\omega'_\lambda,\omega_i\rangle^{-1}e_i,
\quad \omega'_\lambda f_i(\omega'_\lambda)^{-1}=\langle\omega'_\lambda,\omega_i\rangle f_i,
\qquad i\in J,\quad\lambda\in L,
\end{align*}
where $\omega_i=\omega_{\alpha_i}$ and $\omega'_i=\omega'_{\alpha_i}$.
For $L=Q$ and $L=P$, we call it the \emph{root-type} and \emph{weight-type extended Borel algebra}, respectively.
\end{defi}

\begin{prop}\label{aut extendedAn}
  The following morphism $\phi$ gives an nontrivial automorphism of the (root type) extended Borel subalgebra of
   $U_{r, s}(\mathfrak{sl}_{n+1})$.
   \begin{align*}
      & \phi(e_1)=e_n\omega_{n-1}^{-1}\omega_n^{-1}, \quad \phi(e_n)=e_1\omega_1^{-1}\omega_2^{-1}, \\
      & \phi(e_i)=e_{n+1-i}\omega_{n-i}^{-1}\omega_{n+1-i}^{-1}\omega_{n+2-i}^{-1}, \quad 1 < i<n, \\
      & \phi(\omega_j)=(\omega'_{n+1-j})^{-1},   \quad \phi(\omega'_j)=\omega_{n+1-j}^{-1}, \quad 1
      \leq j\leq n.
   \end{align*}
\end{prop}

\begin{proof}
   Recall that  structure constants matrix of $U_{r, s}(\mathfrak{sl}_{n+1})$ is given by
   $$\left(
     \begin{array}{cccccc}
       rs^{-1} & r^{-1} & 1 & \cdots & 1 & 1 \\
       s & rs^{-1} & r^{-1} & \cdots & 1 & 1 \\
       1& s & rs^{-1} & \cdots & 1 & 1 \\
       \vdots & \vdots & \vdots &\vdots  & \vdots &\vdots \\
       1 &1 & 1 & \cdots & rs^{-1} & r^{-1} \\
       1 & 1 & 1 & \cdots & s & rs^{-1} \\
     \end{array}
   \right).$$
   We first verify that $$\phi(\omega_i e_j)=\langle \omega'_j, \omega_i\rangle \phi(e_j\omega_i).$$

 (i) If $i>j+1$ or $i<j-1$,
\[
\phi(\omega_i e_j)=\phi(e_j\omega_i).
\]

(ii) If $i=j+1$,
\begin{align*}
\phi(\omega_{i+1} e_i)
&=\,(\omega'_{n-i})^{-1}e_{n+1-i}\omega_{n-i}^{-1}\omega_{n+1-i}^{-1}\omega_{n+2-i}^{-1} \\
&=\, \langle \omega'_{n-i},\omega_{n+1-i}\rangle e_{n+1-i}\omega_{n-i}^{-1}\omega_{n+1-i}^{-1}\omega_{n+2-i}^{-1}(\omega'_{n-i})^{-1} \\
&=\, \langle \omega'_i,\omega_{i+1}\rangle \phi(e_i\omega_{i+1}) .
\end{align*}

(iii) If $i=j$,
\begin{align*}
\phi(\omega_i e_i)
&=\, (\omega'_{n+1-i})^{-1}e_{n+1-i}\omega_{n-i}^{-1}\omega_{n+1-i}^{-1}\omega_{n+2-i}^{-1} \\
&=\, \langle \omega'_{n+1-i},\omega_{n+1-i}\rangle e_{n+1-i}\omega_{n-i}^{-1}\omega_{n+1-i}^{-1}\omega_{n+2-i}^{-1}(\omega'_{n+1-i})^{-1} \\
&=\, \langle \omega_i',\omega_i\rangle \phi(e_i\omega_i).
\end{align*}

(iv) If $i=j-1$,
\begin{align*}
\phi(\omega_i e_{i+1})
&=\, (\omega'_{n+1-i})^{-1}e_{n-i}\omega_{n-i-1}^{-1}\omega_{n-i}^{-1}\omega_{n-i+1}^{-1} \\
&=\, \langle \omega_{n+1-i}',\omega_{n-i}\rangle e_{n-i}\omega_{n-i-j}^{-1}\omega_{n-i}^{-1}\omega_{n-i+1}^{-1}(\omega'_{n+1-i})^{-1} \\
&=\, \langle \omega_{i+1}',\omega_i\rangle\phi(e_{i+1}\omega_i).
\end{align*}
We can similarly verify that $$\phi(\omega'_i e_j)=\langle \omega'_i, \omega_j\rangle^{-1} \phi(e_j\omega'_i).$$

Then we need to check that \(\phi\) preserves the Serre relations. We may
assume that \(2\leq i\leq n-2\); the cases \(i=1\) and \(i=n\) are checked
similarly. Set
\[
\Omega_i
:=\omega'^{-1}_{n-i}\omega'^{-1}_{n+1-i}\omega'^{-1}_{n+2-i},
\qquad
\Omega_{i+1}
:=\omega'^{-1}_{n-i-1}\omega'^{-1}_{n-i}\omega'^{-1}_{n+1-i}.
\]
Then \(\phi(e_i)=e_{n+1-i}\Omega_i\) and
\(\phi(e_{i+1})=e_{n-i}\Omega_{i+1}\). Moving the group-like factors to the
right, we obtain
\begin{align*}
&\phi\bigl(e_i^2e_{i+1}-(r+s)e_ie_{i+1}e_i
+rs\,e_{i+1}e_i^2\bigr)  \\
&\quad
=s^2\Bigl[
e_{n+1-i}^2e_{n-i}
-(r^{-1}+s^{-1})e_{n+1-i}e_{n-i}e_{n+1-i}
+r^{-1}s^{-1}e_{n-i}e_{n+1-i}^2
\Bigr]\Omega_i^2\Omega_{i+1}  \\
&\quad=0,
\end{align*}
because the expression in brackets is the Serre relation for the adjacent
pair \((n-i,n+1-i)\), multiplied by \(r^{-1}s^{-1}\). Similarly,
\begin{align*}
&\phi\bigl(e_ie_{i+1}^2-(r+s)e_{i+1}e_ie_{i+1}
+rs\,e_{i+1}^2e_i\bigr)  \\
&\quad
=s^2\Bigl[
e_{n+1-i}e_{n-i}^2
-(r^{-1}+s^{-1})e_{n-i}e_{n+1-i}e_{n-i}
+r^{-1}s^{-1}e_{n-i}^2e_{n+1-i}
\Bigr]\Omega_i\Omega_{i+1}^2  \\
&\quad=0.
\end{align*}
From the definition, one can easily verify that $\phi$ is bijective.
\end{proof}

\begin{prop}\label{aut extendedD4}
  The following morphism $\phi$ gives an nontrivial automorphism of the (root type) extended Borel subalgebra of
   $U_{r, s}(\mathfrak{so}_8)$.
   \begin{align*}
     &\phi(e_1)=e_1, \quad \phi(e_2)=e_2, \quad \phi(e_3)=e_4\omega_3\omega_4^{-1}, \quad \phi(e_4)=e_3\omega_3^{-1}\omega_4, \\
     &\phi(\omega_1)=\omega_1, \quad \phi(\omega_2)=\omega_2, \quad \phi(\omega_3)=(\omega_1\omega'_1\omega'_4)^{-1},  \quad \phi(\omega_4)=(\omega_1\omega'_1\omega'_3)^{-1}, \\
     &\phi(\omega'_1)=\omega'_1, \quad \phi(\omega'_2)=\omega'_2, \quad \phi(\omega'_3)=(\omega_1\omega_4\omega'_1)^{-1}, \quad  \phi(\omega'_4)=(\omega_1\omega_3\omega'_1)^{-1}.
   \end{align*}
\end{prop}

\begin{proof}
   Recall that the structure constants matrix of $U_{r, s}(D_4)$ is given by
$$\left(
  \begin{array}{cccc}
    rs^{-1} & r^{-1} & 1 & 1 \\
    s & rs^{-1} & r^{-1} & r^{-1} \\
    1 & s & rs^{-1} & r^{-1}s^{-1} \\
    1 & s & rs & rs^{-1} \\
  \end{array}
\right).$$
   One can easily verify that $\phi$ preserves the defining relations of the extended Borel part.
   For instance,
   \begin{align*}
     \phi(\omega_3 e_1) & =\,(\omega_1\omega'_1\omega'_4)^{-1}e_1  \\
      & =\,e_1(\omega_1\omega'_1\omega'_4)^{-1}  \\
      & =\, \phi(e_1\omega_3).
   \end{align*}
   \begin{align*}
     \phi(\omega_3 e_2) & =\,(\omega_1\omega'_1\omega'_4)^{-1}e_2 \\
      & =\,s^{-1}r^{-1}se_2(\omega_1\omega'_1\omega'_4)^{-1}  \\
      & =\, r^{-1}\phi(e_2\omega_3).
   \end{align*}
      \begin{align*}
     \phi(\omega_3 e_3) & =\,(\omega_1\omega'_1\omega'_4)^{-1}e_4\omega_3\omega_4^{-1} \\
      & =\,rs^{-1}e_4\omega_3\omega_4^{-1}(\omega_1\omega'_1\omega'_4)^{-1}  \\
      & =\, rs^{-1}\phi(e_3\omega_3).
   \end{align*}
   \begin{align*}
     \phi(\omega_3 e_4) & =\, (\omega_1\omega'_1\omega'_4)^{-1}e_3\omega_3^{-1}\omega_4 \\
      & =\, rs e_3\omega_3^{-1}\omega_4(\omega_1\omega'_1\omega'_4)^{-1} \\
      & =\,rs \phi(e_4\omega_3).
   \end{align*}
   \begin{align*}
     \phi(e_3)\phi(e_4) & =\,(e_4\omega_3\omega_4^{-1})(e_3\omega_3^{-1}\omega_4)  \\
      &=\, rs(e_3\omega_3^{-1}\omega_4)(e_4\omega_3\omega_4^{-1})\\
      &=\, rs\phi(e_4)\phi(e_3).
   \end{align*}
   The remaining commutation and quasi-commutation relations follow similarly.
   Next, we show that $\phi$ also preserves Serre relations.
   For instance,
    \begin{align*}
      & \phi\bigl(e_2^2e_4-(r+s)e_2e_4e_2+rse_4e_2^2\bigr)\\
      =\,& e_2^2e_3\omega_3^{-1}\omega_4-(r+s)e_2e_3\omega_3^{-1}\omega_4e_2+rse_3\omega_3^{-1}\omega_4e_2^2  \\
      =\,& \bigl(e_2^2e_3-(r+s)e_2e_3e_2+rse_3e_2^2 \bigr)\omega_3^{-1}\omega_4\\
      =\,&0,
   \end{align*}
       \begin{align*}
      & \phi\bigl(e_2e_4^2-(r+s)e_4e_2e_4+rse_4^2 e_2\bigr)\\
      =\,& e_2e_3\omega_3^{-1}\omega_4e_3\omega_3^{-1}\omega_4-(r+s)e_3\omega_3^{-1}\omega_4 e_2e_3\omega_3^{-1}\omega_4+rse_3\omega_3^{-1}\omega_4e_3\omega_3^{-1}\omega_4e_2  \\
      =\,& r^{-2}\bigl(e_2e_3^2-(r+s)e_3e_2e_3+rse_3^2 e_2 \bigr)\omega_3^{-2}\omega_4^2\\
      =\,&0.
   \end{align*}
   The other Serre relations can also be verified similarly.
From the definition, one can easily verify that $\phi$ is bijective.
\end{proof}

\begin{remark}
We record the action of the automorphism $\phi$ from Proposition \ref{aut extendedD4} on the quantum Lyndon root vectors.
These formulas will be used in Appendix \ref{appA} to derive the commutation relations needed in Appendix \ref{appB}.
   \begin{align*}
    &\phi(\mathcal{E}_{12})=\mathcal{E}_{12}, \quad \phi(\mathcal{E}_{13})=\mathcal{E}_{14}\omega_3\omega_4^{-1}, \quad \phi(\mathcal{E}_{14})=\mathcal{E}_{13}\omega_3^{-1}\omega_4, \\
      &\phi(\mathcal{E}_{23})=\mathcal{E}_{24}\omega_3\omega_4^{-1}, \quad \phi(\mathcal{E}_{24})=\mathcal{E}_{23}\omega_3^{-1}\omega_4, \\
      & \phi(\mathcal{E}_{2\bar{3}})=rs^{-1}\mathcal{E}_{2\bar{3}}, \quad  \phi(\mathcal{E}_{1\bar{3}})=rs^{-1}\mathcal{E}_{1\bar{3}}, \quad \phi(\mathcal{E}_{1\bar{2}})=rs^{-1}\mathcal{E}_{1\bar{2}}.
   \end{align*}
   For instance,
   \begin{align*}
     \phi(\mathcal{E}_{14}) & =\, \phi(\mathcal{E}_{12})\phi(e_4)-s\phi(e_4)\phi(\mathcal{E}_{12}) \\
      & =\,\mathcal{E}_{12}e_3\omega_3^{-1}\omega_4-se_3\omega_3^{-1}\omega_4\mathcal{E}_{12}  \\
      & =\,\bigl(\mathcal{E}_{12}e_3-se_3\mathcal{E}_{12}\bigr)\omega_3^{-1}\omega_4\\
      &=\, \mathcal{E}_{13}\omega_3^{-1}\omega_4,
   \end{align*}
   \begin{align*}
     \phi(\mathcal{E}_{2\bar{3}}) & =\,\phi(\mathcal{E}_{24})\phi(e_3)-r^{-1}\phi(e_3)\phi(\mathcal{E}_{24})  \\
      &=\,\mathcal{E}_{23}\omega_3^{-1}\omega_4e_3\omega_3^{-1}\omega_4-r^{-1}e_3\omega_3^{-1}\omega_4 \mathcal{E}_{23}\omega_3^{-1}\omega_4\\
      &=\,rs^{-1}(r^{-1}s^{-1}\mathcal{E}_{23}e_4-se_4\mathcal{E}_{23}) \\
      &=\, rs^{-1}(r^{-1}s^{-1}\mathcal{E}_{23}e_4-se_4e_2e_3+s^2e_4e_3e_2) \\
     & =\, rs^{-1}\Bigl[r^{-1}s^{-1}\mathcal{E}_{23}e_4-e_2e_4e_3+\mathcal{E}_{24}e_3+r^{-1}e_3(e_2e_4-\mathcal{E}_{24})\Bigr] \\
     & =\, rs^{-1}\Bigl[r^{-1}s^{-1}\mathcal{E}_{23}e_4-r^{-1}s^{-1}e_2e_3e_4+\mathcal{E}_{24}e_3+r^{-1}s^{-1}(e_2e_3-\mathcal{E}_{23})e_4-r^{-1}e_3\mathcal{E}_{24} \Bigr] \\
     & =\,rs^{-1}(\mathcal{E}_{24}e_3-r^{-1}e_3\mathcal{E}_{24})\\
      & =\,rs^{-1}\mathcal{E}_{2\bar{3}},
   \end{align*}
   and
   \begin{align*}
     \phi(\mathcal{E}_{1\bar{3}}) & =\,\phi(e_1)\phi(\mathcal{E}_{2\bar{3}})-s\phi(\mathcal{E}_{2\bar{3}})\phi(e_1) \\
      & =\,rs^{-1}(e_1\mathcal{E}_{2\bar{3}}-s\mathcal{E}_{2\bar{3}}e_1) \\
      & =\,rs^{-1}\mathcal{E}_{1\bar{3}}.
   \end{align*}
   The others can also be derived from definitions directly.
\end{remark}

\iffalse
\begin{remark}\label{rmk anotherLyn}
Let us spell out how a diagram automorphism acts on the convex order in
type $D_4$. Consider the following convex order on $U_{r, s}^+(D_4)$:
\begin{equation}\label{usual convex}
e_1\prec
\mathcal E_{12}\prec
\mathcal E_{13}\prec
\mathcal E_{14}\prec
\mathcal E_{1\bar 3}\prec
\mathcal E_{1\bar 2}\prec
e_2\prec
\mathcal E_{23}\prec
\mathcal E_{24}\prec
\mathcal E_{2\bar 3}\prec
e_3\prec
e_4
\end{equation}
Let $\sigma$ be the nontrivial diagram automorphism interchanging
$\alpha_3$ and $\alpha_4$ and fixing $\alpha_1,\alpha_2$. Since $\sigma$
preserves positive roots and root addition, the transported order is again convex.
That is,
we obtain another convex order (the detailed verification can be seen in the Appendix \ref{appA}) :
\begin{equation}\label{convex order}
e_1
<
\mathcal E_{12}
<
\mathcal E_{14}
<
\mathcal E_{13}
<
\mathcal E_{1\bar 3}
<
\mathcal E_{1\bar 2}
<
e_2
<
\mathcal E_{24}
<
\mathcal E_{23}
<
\mathcal E_{2\bar 3}
<
e_4
<
e_3.
\end{equation}
\end{remark}
\fi

Let $\Gamma_{D_4}$ denote the Dynkin diagram of type $D_4$.
Then the automorphism group of the Dynkin diagram $\operatorname{Aut}(\Gamma_{D_4})\cong S_3$, which admits three generators.
One may naturally ask whether the other two generators induce automorphisms of the (root type) extended Borel subalgebra of $U_{r,s}(\mathfrak{so}_8)$.
The following proposition provides a negative answer to this question.
\begin{prop}\label{prop aut D4}
Let $\phi$ and $\psi$ be the two nontrivial generators of $\operatorname{Aut}(\Gamma_{D_4})$ defined by
\[
\phi(\alpha_1)=\alpha_3,\quad \phi(\alpha_2)=\alpha_2,\quad \phi(\alpha_3)=\alpha_1,\quad \phi(\alpha_4)=\alpha_4,
\]
\[
\psi(\alpha_1)=\alpha_4,\quad \psi(\alpha_2)=\alpha_2,\quad \psi(\alpha_3)=\alpha_3,\quad  \psi(\alpha_4)=\alpha_1.
\]
Then neither $\phi$ nor $\psi$ lifts to an automorphism of the root-type extended Borel algebra of $U_{r,s}(\mathfrak{so}_8)$.
\end{prop}

\begin{proof}
If $\phi$ lifts to an automorphism of the algebra, we retain the same symbol $\phi$ for the lifted automorphism.
We first determine the images of the group-like elements: $\phi(\omega_i)$ and $\phi(\omega'_i)$.
The images of the simple root vectors $\phi(e_i)$ are determined up to multiplication by nonzero scalars and group-like elements.
Since such factors do not affect the conjugation relations determining the images of the group-like elements,
we normalize the images and assume
$$\phi(e_1)=e_3,\quad \phi(e_2)=e_2,\quad \phi(e_3)=e_1,\quad \phi(e_4)=e_4.$$
(This is only a temporary assumption.
The precise formula of these $\phi(e_i)$ will be given in Theorem \ref{prop aut D4weight}).
We take the determination of $\phi(\omega_3)$ as an example.

According to \cite{Fleury 97}, every invertible element is a nonzero scalar multiple of a Laurent monomial in the group-like elements.
Without loss of generality, we assume that $\phi(\omega_3)=\omega_1^{a_1}\omega_2^{a_2}\omega_3^{a_3}\omega_4^{a_4}(\omega'_1)^{b_1}(\omega'_2)^{b_2}(\omega'_3)^{b_3}(\omega'_4)^{b_4}$.
From the definition, we know that
\begin{align*}
 \phi(\omega_3e_1)=\phi(e_1\omega_3),\qquad & \phi(\omega_3e_2)=r^{-1}\phi(e_2\omega_3), \\
  \phi(\omega_3e_3)=rs^{-1}\phi(e_3\omega_3), \qquad & \phi(\omega_3e_4)=rs\phi(e_4\omega_3).
\end{align*}
Then we have
\begin{align*}
   & s^{a_2}(rs^{-1})^{a_3}(r^{-1}s^{-1})^{a_4}r^{b_2}(r^{-1}s)^{b_3}(r^{-1}s^{-1})^{b_4}=1,  \\
   & s^{a_1}(rs^{-1})^{a_2}(r^{-1})^{a_3}(r^{-1})^{a_4}r^{b_1}(r^{-1}s)^{b_2}(s^{-1})^{b_3}(s^{-1})^{b_4}=r^{-1},\\
   & (rs^{-1})^{a_1}(r^{-1})^{a_2}(r^{-1}s)^{b_1}(s^{-1})^{b_2}=rs^{-1}, \\
   & s^{a_2}(rs)^{a_3}(rs^{-1})^{a_4}r^{b_2}(rs)^{b_3}(r^{-1}s)^{b_4}=rs.
\end{align*}
Comparing the corresponding coefficients of $r$ and $s$, we obtain a system consisting of eight equations in eight variables.
The coefficient matrix is invertible, so we can obtain an unique solution:
$\phi(\omega_3)=\omega_1^{\frac{1}{2}}\omega_3^{\frac{1}{2}}(\omega'_1)^{-\frac{1}{2}}(\omega'_3)^{\frac{1}{2}}.$
Similarly, we can uniquely determine the images of all group-like elements, up to some scalars:
   \begin{align*}
     &\phi(\omega_1)=\omega_1^{\frac{1}{2}}\omega_3^{\frac{1}{2}}(\omega'_1)^{\frac{1}{2}}(\omega'_3)^{-\frac{1}{2}}, \quad \phi(\omega_2)=\omega_1^{\frac{1}{2}}\omega_3^{\frac{1}{2}}(\omega'_1)^{\frac{1}{2}}(\omega'_2)^{-1}(\omega'_3)^{\frac{1}{2}}, \\
&\phi(\omega_3)=\omega_1^{\frac{1}{2}}\omega_3^{\frac{1}{2}}(\omega'_1)^{-\frac{1}{2}}(\omega'_3)^{\frac{1}{2}},  \quad \phi(\omega_4)=\omega_1^{-\frac{1}{2}}\omega_2\omega_3^{-\frac{1}{2}}\omega_4(\omega'_1)^{-\frac{1}{2}}\omega'_2(\omega'_3)^{-\frac{1}{2}}, \\
     &\phi(\omega'_1)=\omega_1^{\frac{1}{2}}\omega_3^{-\frac{1}{2}}(\omega'_1)^{\frac{1}{2}}(\omega'_3)^{\frac{1}{2}}, \quad \phi(\omega'_2)=\omega_1^{\frac{1}{2}}\omega_2^{-1}\omega_3^{\frac{1}{2}}(\omega'_1)^{\frac{1}{2}}(\omega'_3)^{\frac{1}{2}}, \\
&\phi(\omega'_3)=\omega_1^{-\frac{1}{2}}\omega_3^{\frac{1}{2}}(\omega'_1)^{\frac{1}{2}}(\omega'_3)^{\frac{1}{2}}, \quad  \phi(\omega'_4)=\omega_1^{-\frac{1}{2}}\omega_2\omega_3^{-\frac{1}{2}}(\omega'_1)^{-\frac{1}{2}}\omega'_2(\omega'_3)^{-\frac{1}{2}}\omega'_4.
   \end{align*}
$\phi$ can be viewed as the morphism of the weight lattice $P\oplus P$.
We will take  $\phi(\omega_1)=\omega_1^{\frac{1}{2}}\omega_3^{\frac{1}{2}}(\omega'_1)^{\frac{1}{2}}(\omega'_3)^{-\frac{1}{2}}$ as an illustration.

As is well-known, for type $D_4$,
the simple roots $\alpha_i\, (1\leq i\leq 4)$ and the fundamental weights $\lambda_i\, (1\leq i\leq 4)$ satisfy the following relations:
\begin{align}\label{equ rootweight}
\alpha_1&=2\lambda_1-\lambda_2, \nonumber\\
\alpha_2&=-\lambda_1+2\lambda_2-\lambda_3-\lambda_4,\\
\alpha_3&=-\lambda_2+2\lambda_3, \nonumber\\
\alpha_4&=-\lambda_2+2\lambda_4. \nonumber
\end{align}
From (\ref{equ rootweight}), we have
\begin{align*}
   \tfrac12\alpha_1+\tfrac12\alpha_3&=\,\tfrac12(2\lambda_1-\lambda_2)+\tfrac12(-\lambda_2+2\lambda_3) \\
   & =\,\lambda_1-\lambda_2+\lambda_3\in P.\\
   \tfrac12\alpha_1-\tfrac12\alpha_3&=\,\tfrac12(2\lambda_1-\lambda_2)-\tfrac12(-\lambda_2+2\lambda_3) \\
   & =\,\lambda_1-\lambda_3\in P.
\end{align*}
The other formulas can be checked similarly.
Since the corresponding exponent matrices with respect to the basis
$\{\omega_1,\cdots,\omega_4,\omega'_1,\cdots,\omega'_4\}$
belong to $GL(8,\mathbb{Z})$,
they determine an automorphism of $P\oplus P$.
However, these images do not belong to the root lattice $Q\oplus Q$, since half-integral exponents occur.
Therefore, no automorphism of the root-type extended Borel algebra can induce these transformations on the group-like elements.
If $\psi$ lifts to an automorphism of the algebra, we retain the same symbol $\psi$ for the lifted automorphism.
Similarly, we can uniquely determine the images of all group-like elements, up to some scalars:
  \begin{align*}
     &\psi(\omega_1)=\omega_1^{\frac{1}{2}}\omega_4^{\frac{1}{2}}(\omega'_1)^{\frac{1}{2}}(\omega'_4)^{-\frac{1}{2}}, \quad \psi(\omega_2)=\omega_1^{-\frac{1}{2}}\omega_4^{-\frac{1}{2}}(\omega'_1)^{-\frac{1}{2}}(\omega'_2)^{-1}(\omega'_4)^{-\frac{1}{2}}, \\
&\psi(\omega_3)=\omega_1^{-\frac{1}{2}}\omega_2^{-1}\omega_3\omega_4^{-\frac{1}{2}}(\omega'_1)^{-\frac{1}{2}}(\omega'_2)^{-1}(\omega'_4)^{-\frac{1}{2}},  \quad \psi(\omega_4)=\omega_1^{\frac{1}{2}}\omega_4^{\frac{1}{2}}(\omega'_1)^{-\frac{1}{2}}(\omega'_4)^{\frac{1}{2}}, \\
     &\psi(\omega'_1)=\omega_1^{\frac{1}{2}}\omega_4^{-\frac{1}{2}}(\omega'_1)^{\frac{1}{2}}(\omega'_4)^{\frac{1}{2}}, \quad \psi(\omega'_2)=\omega_1^{-\frac{1}{2}}\omega_2^{-1}\omega_4^{-\frac{1}{2}}(\omega'_1)^{-\frac{1}{2}}(\omega'_4)^{-\frac{1}{2}}, \\
&\psi(\omega'_3)=\omega_1^{-\frac{1}{2}}\omega_2^{-1}\omega_4^{-\frac{1}{2}}(\omega'_1)^{-\frac{1}{2}}(\omega'_2)^{-1}\omega'_3(\omega'_4)^{-\frac{1}{2}},\quad  \psi(\omega'_4)=\omega_1^{-\frac{1}{2}}\omega_4^{\frac{1}{2}}(\omega'_1)^{\frac{1}{2}}(\omega'_4)^{\frac{1}{2}}.
   \end{align*}
One can easily verify that $\psi$ is an automorphism of $P\oplus P$,
and similarly, $\psi$ cannot be extended to an automorphism of the root-type extended Borel algebra.
\end{proof}

\begin{theorem}\label{prop aut D4weight}
Let $\sigma\in\operatorname{Aut}(\Gamma_{D_4})$ be a Dynkin diagram automorphism.
Then $\sigma$ lifts to an automorphism of the extended Borel algebra
if and only if the extended Borel algebra is of weight-type.

Equivalently, $\operatorname{Aut}(\Gamma_{D_4})$ is fully realized on the weight-type extended Borel algebra,
whereas in the root-type case only the subgroup generated by the Dynkin diagram automorphism interchanging $\alpha_3$ and $\alpha_4$ is realized.
\end{theorem}

\begin{proof}
Define two algebra endomorphisms $\phi$ and $\psi$ of the weight-type extended Borel algebra as follows:
The images of the positive simple root vectors are given by
\begin{align*}
   & \phi(e_1)=e_3\omega_1\omega_3^{-1},\quad \phi(e_2)=e_2,\quad \phi(e_3)=e_1\omega_1^{-1}\omega_3,\quad \phi(e_4)=e_4; \\
   & \psi(e_1)=e_4\omega_1\omega_4^{-1},\quad  \psi(e_2)=e_2,\quad  \psi(e_3)=e_3,\quad \psi(e_4)=e_1\omega_1^{-1}\omega_4.
\end{align*}
The images on the group-like elements are given in the proof of Proposition \ref{prop aut D4}.
We check that $\phi$ indeed preserves the defining relations of the $U_{r, s}^+(\mathfrak{so}_8)$.
First we verify the commutation relations among the positive generators.
\begin{align*}
  \phi(e_3e_4)\, &=e_1\omega_1^{-1}\omega_3e_4 \\
   &=rs e_4e_1\omega_1^{-1}\omega_3 \\
  &=rs\phi(e_4e_3).\\
 \phi(e_1e_3)\, &=e_3\omega_1\omega_3^{-1}e_1\omega_1^{-1}\omega_3 \\
   &=rs^{-1}r^{-1}s(e_1\omega_1^{-1}\omega_3)(e_3\omega_1\omega_3^{-1}) \\
  &=\phi(e_3e_1).\\
 \phi(e_1e_4)\, &=e_3\omega_1\omega_3^{-1}e_4 \\
   &=e_4e_3\omega_1\omega_3^{-1} \\
  &=\phi(e_4e_1).
\end{align*}
The defining relations involving the group-like generators $\phi(\omega_i)$
and $\phi(\omega_i')$ have already been verified in Proposition~\ref{prop aut D4}.
It remains to verify the Serre relations.
\begin{align*}
   & \phi\bigl(e_1^2e_2-(r+s)e_1e_2e_1+rse_2e_1^2\bigr) \\
  =\, & e_3\omega_1\omega_3^{-1}e_3\omega_1\omega_3^{-1}e_2-(r+s)e_3\omega_1\omega_3^{-1}e_2e_3\omega_1\omega_3^{-1}+rse_2e_3\omega_1\omega_3^{-1}e_3\omega_1\omega_3^{-1} \\
  =\, &s^2\bigl[e_2e_3^2-(r+s)e_3e_2e_3+rse_3^2e_2\bigr]\omega_1^{-2}\omega_3^{-2} \\
  =\, & 0.\\
   & \phi\bigl(e_1e_2^2-(r+s)e_2e_1e_2+rse_2^2 e_1\bigr) \\
  =\, & e_3\omega_1\omega_3^{-1}e_2^2-(r+s)e_2e_3\omega_1\omega_3^{-1}e_2+rse_2^2e_3\omega_1\omega_3^{-1} \\
  =\, &r^2s^2\bigl[e_3e_2^2-(r^{-1}+s^{-1})e_2e_3e_2+r^{-1}s^{-1}e_2^2e_3\bigr]\omega_1\omega_3^{-1} \\
  =\, & 0.\\
& \phi\bigl(e_2^2e_3-(r+s)e_2e_3e_2+rse_3e_2^2\bigr) \\
  =\, & e_2^2e_1\omega_1^{-1}\omega_3-(r+s)e_2e_1\omega_1^{-1}\omega_3e_2+rse_1\omega_1^{-1}\omega_3e_2^2 \\
  =\, &r^{-1}s^{-1}\bigl[e_1e_2^2-(r+s)e_2e_1e_2+rse_2^2e_1\bigr]\omega_1^{-1}\omega_3 \\
  =\, & 0.\\
& \phi\bigl(e_2e_3^2-(r+s)e_3e_2e_3+rse_3^2 e_2\bigr) \\
  =\, & e_2e_1\omega_1^{-1}\omega_3e_1\omega_1^{-1}\omega_3-(r+s)e_1\omega_1^{-1}\omega_3e_2e_1\omega_1^{-1}\omega_3+rse_1\omega_1^{-1}\omega_3e_1\omega_1^{-1}\omega_3e_2 \\
  =\, &r^{-2}\bigl[e_1^2 e_2-(r+s)e_1e_2e_1+rse_2e_1^2\bigr]\omega_1^{-2}\omega_3^2 \\
  =\, & 0.
\end{align*}
Hence we finish the verification for $\phi$.
One can easily check that $\psi$ also preserve the defining relations of the $U_{r, s}^+(\mathfrak{so}_8)$ analogously.
Moreover, both $\phi$ and $\psi$ are bijective, and hence are automorphisms of the weight-type extended Borel algebra.

Together with Proposition~\ref{aut extendedD4}, this shows that every element of $\operatorname{Aut}(\Gamma_{D_4})\cong S_3$
lifts to an automorphism of the weight-type extended Borel algebra,
whereas only the subgroup generated by the Dynkin diagram automorphism interchanging $\alpha_3$ and $\alpha_4$ does so in the root-type extended Borel algebra.
\end{proof}

\section{Representation-theoretic distinctions stemming from commutative and quasi-commutative relations}\label{sec repdis}

We will also reveal the algebraic non-isomorphic relationship between two-parameter and one-parameter counterparts caused by the ``{\it type $D$ phenomenon}",
which reflects the fact that {\it some commutative relations in $U_q(\mathfrak{g})$ have to turn into quasi-commutative ones in $U_{r,s}(\mathfrak{g})$ with suitable powers of factor $rs$ in two-parameter setup}.
Indeed, this phenomenon not only appears in the definition  of $U_{r,s}\mathcal(\mathfrak{so}_{2n})$ (see \cite{BGH 2006}),
but also extensively in two-parameter quantum affine setting,
see their structure constant matrices except for $A_n^{(1)}$ (\cite{Zhang phd 2007, HuZhang2014}),
such as $U_{r,s}\mathcal(\widehat {\mathfrak{so}_{N}})$ with $N\in\{2n, 2n+1\}$ $\, U_{r,s}\mathcal(\widehat {\frak{sp}_{2n}}), \, U_{r,s}\mathcal(\widehat {E_N})$ with $N=6,7,8$, $ U_{r,s}\mathcal(\widehat {F_4})$, $ U_{r,s}\mathcal(\widehat {G_2})$, etc (See \cite{BGH 2006, HRZ 2008, HuZhang2014, Zhang phd 2007}).
Particularly, one should see Axioms $(D1), (D6_1), (D6_2)$ of Drinfeld realization due to \cite{HRZ 2008} (consult more remarks in \cite{HuZhang2014}).
We refer to it as the “{\it type $D$-phenomenon}”,
since this firstly appeared in our definition of $U_{r,s}\mathcal(\frak{so}_{2n})$ (see \cite{BGH 2006}).
We explain it explicitly by showing the Dynkin diagram of type $D$:

\vspace{1em}
\hspace{5em}
$$\begin{tikzpicture}[scale=0.6]
\draw (-2,0)--(0,0);
\draw (0,0)-- (2,0);
\draw (2,0)-- (2.7,0);
\node at (3,0) {-};
\node at (3.4,0) {-};
\node at (3.8,0) {-};
\node at (4.15,0) {-};
\draw (4.3,0)-- (5,0);
\draw (5,0) -- (7,-1.5);
\draw (5,0) -- (7,1.5);
\draw[fill=white] (-2 cm, 0 cm) circle (.15cm) node[below=3pt]{$1$};
\draw[fill=white] (0 cm, 0 cm) circle (.15cm) node[below=3pt]{$2$};
\draw[fill=white] (2 cm, 0 cm) circle (.15cm) node[below=3pt]{$3$};
\draw[fill=white] (5 cm, 0 cm) circle (.15cm) node[below=5pt]{$n-2$};
\draw[fill=white] (7 cm, 1.5 cm) circle (.15cm) node[above=3pt]{$n-1$};
\draw[fill=white] (7 cm, -1.5 cm) circle (.15cm) node[below=6pt]{$n$};
\end{tikzpicture}$$

It should be noted  that the dots labeled $n-1$ and $n$ are disconnected.
In other words, their corresponding simple root vectors are commutative in $U_q(\mathfrak{so}_{2n})$:
$e_{n-1}e_n=e_ne_{n-1}$.
However, in the two-parameter case \cite{BGH 2006}, this relation turns to $e_{n-1}e_n=rs\,e_ne_{n-1}$.
Obviously the former relation is the degenerate case of the latter when $r=q=s^{-1}$, i.e., $rs=1$.
Nevertheless, $e_{n-1}$ and $e_n$ are no longer commutative in $U_{r, s}(\mathfrak{so}_{2n})$.
This difference in the defining relations motivates a comparison of their representation theory beyond the category of type~1 weight modules.
Namely, the categories of weight modules cannot detect such type $D$-phenomenon (especially occurring in the two-parameter quantum affine setup). In this regard, we shall point out that it will be very important topics in their representation theory for $U_{r,s}(\widehat{\mathfrak g})$ if one starts and expands to study the Whittaker categories.

\begin{prop}\cite{HuPeiCA 2012}
  Set $q^2=rs^{-1}$. As braided tensor categories, we have the following equivalence:
   $$\mathcal O^{r, s}\cong \mathcal O^{q, q^{-1}},$$
   where $\mathcal O^{r, s}$ and $ \mathcal O^{q, q^{-1}}$ represent the type $1$ weight module category of $U_{r, s}(\frak{g})$ and $U_{q, q^{-1}}(\frak{g})$, respectively.
\end{prop}

Consider the subalgebra $\mathcal{B}_{rs}=\langle e_{n-1}, e_n\rangle \subseteq U_{r, s}(\frak{so}_{2n})$
and $\mathcal{B}_1=\langle e_{n-1}, e_n\rangle \subseteq U_{q, q^{-1}}(\frak{so}_{2n})$.
The next Proposition will reveal the difference between their module categories.

\begin{prop}\label{prop typeD1}
  Define $\mathcal{B}_p:=\mathbb{C}\langle x,y\rangle/(xy-pyx)$.

  $(1)$ If $p\in \mathbb{C}^\times$ is not a root of unity, there exists an infinite-dimensional faithful simple module over $\mathcal{B}_p$.

  $(2)$ If $p\in \mathbb{C}^\times$ is an $\ell$-th primitive root of unity, there exists an $\ell$-dimensional simple module over $\mathcal{B}_p$.
\end{prop}

\begin{proof}
  (1)  Assume $V$ is an infinite-dimensional vector space with a basis $\{v_i \mid i\in \mathbb{Z}\}$.
  We define the actions of $x$ and $y$ on the basis vectors as follows:
  $$x\cdot v_i=v_{i+1}, \quad y\cdot v_i=p^{-i}v_{i-1}.$$
  Obviously, $V$ is an infinite-dimensional $\mathcal{B}_p$-module.
  Let $M\subset V$ be a nonzero $\mathcal{B}_p$-submodule and take a nonzero element $m=\sum_{i=j}^{k}a_iv_i\in M$.
  For any non-negative integer $t$, we have
  \begin{align*}
    (xy)^t\cdot m=\sum_{i=j}^{k}a_i(xy)^t\cdot v_i=\sum_{i=j}^{k} p^{-it}a_iv_i.
  \end{align*}
  Thus we can get all non-zero summands $a_i v_i\; (j\leq i\leq k)$ in $M$ (since $p$ is not a root of unity).
  By the action of $x$ and $y$, we thus get $v_i\in M$ for $i\in \mathbb{Z}$.
  Therefore $V=M$ is simple. Similarly, we can also prove this module is faithful.

  (2) Let $V$ be an $\ell$-dimensional vector space with a basis $\{v_i \mid 0\leq i\leq \ell-1\}$.
  We define
  $$x\cdot v_i=v_{i+1}\; \text{for}\; 0\leq i\leq \ell-2;  \quad x\cdot v_{\ell-1}=v_0,$$
  $$y\cdot v_i=p^{-i}v_{i-1} \;\text{for}\; 1\leq i\leq \ell-1; \quad y\cdot v_0=(-1)^{\ell+1}v_{\ell-1}.$$
  Then we can prove $V$ is a simple module similarly.
\end{proof}

\begin{remark}
Assume $p=1$, and let $V = \{v_i \mid i\in \mathbb{Z}\}$ be a $\mathcal{B}_1$-module with the action
$$x\cdot v_i = av_{i+1},\qquad y\cdot v_i = bv_{i-1}, \quad a, b\in \mathbb{C}.$$
Then $V$ is neither simple nor Artinian.
\end{remark}

\begin{proof}
  For each non-negative integer $k$, define
$$V_k=\operatorname{Span}\{\,v_i-v_{i+2^k}\mid i\in\mathbb Z\,\}.$$
Then $V_0$ is obviously a  $\mathcal{B}_1$-submodule of $V$, and each $V_k$  is a $\mathcal{B}_1$-submodule of $V_{k-1}$.
Thus we have the (infinite) strictly decreasing chain
\[
V\supsetneq V_0\supsetneq V_1\supsetneq V_2\supsetneq \cdots \supsetneq V_n \supsetneq \cdots.
\]
Consequently, $V$ is neither simple nor Artinian as a $\mathcal{B}_1$-module.
\end{proof}

\begin{remark}
The infinite-dimensional simple module in Proposition \ref{prop typeD1}\,(1) can be recovered from an induced module after localization.
Let \(S_x=\{x^i\mid i\geq 0\}\).
Since \(x\) is a normal element of \(\mathcal B_p\), the set \(S_x\) is a two-sided Ore set.
Hence we may consider the Ore localization \(A=\mathcal B_p[S_x^{-1}]\).
Since \(x\) acts bijectively on \(V\), the \(\mathcal B_p\)-module \(V\) naturally extends to an \(A\)-module.

Set \(X=x\) and \(Y=xy\).
Then \(A\simeq \mathbb C\langle X^{\pm1},Y\rangle/(XY-pYX)\).
Moreover, \(Yv_k=p^{-k}v_k\) for every \(k\in \mathbb Z\).
Let \(\mathbb C_k=\mathbb C v_k\) be the one-dimensional \(\mathbb C[Y]\)-module on which \(Y\) acts by the scalar \(p^{-k}\).
Define
\[
\Phi_k:A\otimes_{\mathbb C[Y]}\mathbb C_k\longrightarrow V,\qquad
X^iY^j\otimes v_k\longmapsto p^{-jk}v_{i+k},
\]
where \(i\in\mathbb Z\) and \(j\geq 0\).
This map is well-defined and is an \(A\)-module isomorphism.
Indeed, the tensor relation gives \(X^iY^j\otimes v_k=X^i\otimes Y^jv_k=p^{-jk}X^i\otimes v_k\),
and the inverse is given by \(v_{i+k}\mapsto X^i\otimes v_k\).

Thus the module \(V\) is obtained by inducing the one-dimensional \(\mathbb C[Y]\)-module \(\mathbb C_k\) to the localization \(A\).
If a different eigenvector \(v_l\) is chosen, the corresponding induced module is again naturally isomorphic to \(V\).
Hence the realization is independent of the choice of \(v_k\) up to this natural isomorphism.
\end{remark}

\begin{theorem}\label{thm type D isoclass}
  Assume that $\mathcal{B}_p=\mathbb{C}\langle x,y\rangle/(xy-pyx)$.

  $(1)$ If $p\in\mathbb{C}^\times$ is not a root of unity,
  then all simple $\mathcal{B}_p$-modules are either one-dimensional or infinite-dimensional.

  $(2)$ If $p\in \mathbb{C}^\times$ is an $\ell$-th primitive root of unity,
  then all simple $\mathcal{B}_p$-modules are either one-dimensional or $\ell$-dimensional.
  Moreover, $\varphi: \mathcal{B}_p\rightarrow M_\ell(\mathbb{C})$ is surjective.

  $(3)$ If $p=1$, every simple module over $\mathcal{B}_1$ has dimension one.
\end{theorem}

\begin{proof}
  (1) To classify all the simple modules over $\mathcal{B}_p$, we investigate their annihilators, which are exactly primitive ideals of $\mathcal{B}_p$.
  According to Example II 7.2 in \cite{Brown}, all the primitive ideals of $\mathcal{B}_p$ are $\langle 0\rangle, \langle x, y-\beta\rangle, \langle x-\alpha, y\rangle$, and $\langle x, y\rangle $  (where $\alpha, \beta\in \mathbb{C}^\times)$.
  If the annihilator is $\langle x-\alpha, y\rangle$ or $\langle x, y-\beta\rangle$, then the corresponding simple module is $1$-dimensional (That is, $x$ or $y$ acts as non-zero scalar multiplication, while another indeterminate acts as zero).

  If the annihilator is $\langle 0 \rangle$, we claim that these faithful simple modules must be infinite-dimensional.
  If not, suppose that $V$ is $n$-dimensional simple faithful module.
Hence the actions of $x$ and $y$ on $V$ are represented by matrices $ X,Y\in M_n(\mathbb{C})$ satisfying $XY=pYX$.
Therefore, $ (1-p^n)\det(X)\det(Y)=0$.
Since $p$ is not a root of unity, $1-p^n\neq0$, we have $\det(X)\det(Y)=0$.
Thus at least one of $X$ and $Y$ is singular. Without loss of generality, assume that $\det(X)=0$.
Then $\ker X\neq0$.
For any $v\in\ker X$, we have $X(Yv) = pY(Xv) = 0$, which implies that $ Y(\ker X)\subseteq\ker X$.
Since $X$ annihilates $\ker X$, it follows that $\ker X$ is a nonzero $\mathcal{B}_p$-submodule of $V$.
By the simplicity of $V$, we obtain $\ker X=V$, and hence $ X=0$.
Consequently, every finite-dimensional simple $\mathcal{B}_p$-module factors through either $\mathcal{B}_p/(x)\cong\mathbb{C}[y]$.
 Since every simple module over $\mathbb{C}[y]$ is one-dimensional, every finite-dimensional simple $\mathcal{B}_p$-module is one-dimensional and therefore cannot be faithful.
Hence a faithful simple $\mathcal{B}_p$-module must be infinite-dimensional.
  One concrete construction is given in Proposition \ref{prop typeD1}.

  (2) Since $x^\ell$ and $y^\ell$ are central elements of $\mathcal{B}_p$, they act as scalar multiplication on simple modules.
  Assume that $x^\ell= \alpha \cdot\text{id}$, and $y^\ell=\beta \cdot \text{id}$, where $\alpha, \beta\in \mathbb{C}$.
  If $\alpha\beta=0$, then it is corresponding to $1$-dimensional simple module.
  If $\alpha\beta\neq 0$, we may assume $\alpha=\beta=1$ without loss of generality.
  Then we have the following isomorphism:
  $$\mathcal{B}_p/\langle x^\ell-1, y^\ell-1\rangle\cong M_\ell(\mathbb{C}).$$
  Indeed, define $\varphi: \mathcal{B}_p\rightarrow   M_\ell(\mathbb{C})$ by
  $$\varphi(x):=\left(
                 \begin{array}{ccccc}
                   0 & 0 & \cdots  & 0 & 1 \\
                   1 & 0 & \cdots  & 0 & 0 \\
                   0 & 1 & \cdots  & 0 & 0 \\
                   \vdots & \vdots &    & \vdots & \vdots \\
                   0 & 0 & \cdots  & 1 & 0 \\
                 \end{array}
               \right), \quad
  \varphi(y):=\left(
                \begin{array}{ccccc}
                  0 & p^{-1} & 0 & \cdots & 0 \\
                  0 & 0 & p^{-2} & \cdots & 0 \\
                  \vdots & \vdots & \vdots &  & \vdots \\
                  0 & 0 & 0 & \cdots & p^{-(\ell-1)} \\
                  (-1)^{\ell-1} & 0 & 0 & \cdots & 0 \\
                \end{array}
              \right).
               $$
  Obviously, $\varphi$ is a homomorphism and $\text{Ker}\varphi=\langle x^\ell-1, y^\ell-1\rangle$.
  From Proposition \ref{prop typeD1},
  we know that  $S:=\langle \varphi(x), \varphi(y)\rangle$ is a subalgebra of $T:=\text{End}_\mathbb{C}(\mathbb{C}^\ell)=M_\ell(\mathbb{C})$, and $\mathbb{C}^\ell$ is a simple $S$-module.
  Using the Double Centralizer Theorem, we conclude that
  $$T=C_T(C_T(S))=S.$$
  Thus, $\varphi$ is surjective.
  In this case, all the simple modules over $\mathcal{B}_p/\langle x^\ell-1, y^\ell-1\rangle $ have dimension $\ell$.
  One concrete construction is given in Proposition \ref{prop typeD1}.

  (3) Every simple module over $\mathcal{B}_1\cong \mathbb{C}[x,y]$ has the form of $\mathbb{C}[x, y]/I$,
  where $I$ is a maximal ideal of $\mathbb{C}[x,y]$.
  By Hilbert's Nullstellensatz, $I=(x-a, y-b), a, b\in \mathbb{C}$.
  That is to say, every simple module over $\mathcal{B}_1$ has dimension one.
\end{proof}

The preceding comparison of simple modules shows that the representation theories of $\mathcal B_p$ and $\mathcal B_1$ behave quite differently.
However, the vector-space dimensions of simple modules are not, by themselves, invariants of an equivalence of module categories.
To overcome this obstruction, we consider their categorical centers.
Recall that, for a $\mathbb C$-linear category $\mathcal C$,
its categorical center $\mathfrak Z(\mathcal C):=\operatorname{End}(\operatorname{Id}_{\mathcal C})$,
the algebra of natural endomorphisms of the identity functor \cite{Khal09}.
This categorical center is invariant under equivalence of categories.
Moreover, for a unital $\mathbb C$-algebra $A$, one has a natural isomorphism $\mathfrak Z(A\text{-Mod})\cong Z(A)$,
where $Z(A)$ denotes the ordinary center of $A$.
Therefore, if $\mathcal B_p\text{-Mod}$ and $\mathcal B_1\text{-Mod}$ were equivalent,
then their algebra centers would be isomorphic.

\begin{prop}\label{prop typeD2}
Assume that $p$ is not a root of unity.
Then $\mathcal B_p\text{-}\mathrm{Mod}$ and $\mathcal B_1\text{-}\mathrm{Mod}$ are not equivalent as abelian categories.
\end{prop}

\begin{proof}
It suffices to prove $Z(\mathcal B_p)\ncong Z(\mathcal{B}_1)$.
Now \(\mathcal B_1=\mathbb C[x,y]\), so \(Z(\mathcal B_1)=\mathbb C[x,y]\).
On the other hand, since \(\mathcal B_p=\mathbb C\langle x,y\rangle/(xy-pyx)\) and \(p\) is not a root
of unity, we have \(Z(\mathcal B_p)=\mathbb C\).
Indeed, writing an element \(z\in \mathcal B_p\) as \(z=\sum c_{ij}x^iy^j\),
the equality \(xz=zx\) forces \(j=0\) whenever \(c_{ij}\neq 0\), and the equality \(yz=zy\) then forces \(i=0\) whenever \(c_{ij}\neq 0\).
Thus \(z\) is a scalar.

Therefore \(Z(\mathcal B_p)\simeq \mathbb C\), whereas \(Z(\mathcal B_1)\simeq \mathbb C[x,y]\).
Hence \(\mathcal B_p\text{-}\mathrm{Mod}\) and \(\mathcal B_1\text{-}\mathrm{Mod}\) are not equivalent as abelian categories.
\end{proof}

\section{Induced modules from the distinguished subalgebra to $U_{r,s}(\mathfrak{so}_{2n})$}
\label{sec induced}

In this section, we investigate the whole two-parameter quantum group $U_{r,s}(\mathfrak{so}_{2n})$ through induced modules.
The starting point is the classification of simple modules over the distinguished subalgebra $\mathcal B_{rs}$,
generated by $e_{n-1}$ and $e_n$, obtained in Theorem \ref{thm type D isoclass}.
Using Frobenius reciprocity, we prove that induced modules coming from one-dimensional simple $\mathcal B_{rs}$-modules admit no nonzero
homomorphisms to those coming from infinite-dimensional simple $\mathcal B_{rs}$-modules.
This provides a representation-theoretic distinction between the two-parameter and one-parameter cases.

\begin{prop}[Frobenius reciprocity]\label{prop Frobenius}
Assume that $\mathcal{B}$ is a subalgebra of $\mathcal{A}$.
For any $\mathcal{A}$-module $M$ and $\mathcal{B}$-module $N$, we have the following canonical isomorphism:
$$\text{Hom}_\mathcal{A}(\mathcal{A}\otimes_{\mathcal{B}} N, M)\cong \text{Hom}_{\mathcal{B}}(N, \text{Res}_{\mathcal{B}}^{\mathcal{A}} M).$$
\end{prop}

\begin{defi}\label{def Ars}
Consider the algebra $\mathcal{A}_{r,s} :=U_{r,s}(\mathfrak{so}_{2n})= \langle e_i, f_i, \omega_i^{\pm 1}, (\omega'_i)^{\pm 1}\mid 1\leq i\leq n\rangle$, which satisfies Definition \ref{def DJ}.
Its structure constants matrix is given by
$$\left(
  \begin{array}{cccccc}
    rs^{-1} & r^{-1} & 1 & \cdots & 1 & 1 \\
    s & rs^{-1} & r^{-1} & \cdots & 1 & 1 \\
    1 & s & rs^{-1} & \cdots & 1 & 1 \\
    \vdots & \vdots & \vdots & \vdots & \vdots & \vdots \\
    1 & 1 & 1 & \cdots &  rs^{-1} &  r^{-1}s^{-1} \\
    1 & 1 & 1 & \cdots & rs & rs^{-1} \\
  \end{array}
\right)$$
\end{defi}
In particular, let $r=q=s^{-1}$, we obtain the one-parameter quantum group $\mathcal{A}_1=U_{q, q^{-1}}(\mathfrak{so}_{2n})$.

Recall that given a Hopf algebra $H$ and a Hopf $2$-cocycle $\sigma$, one can construct a new Hopf algebra $H^\sigma$ \cite{Doi 1994}, with the same comultiplication, while the new multiplication can be defined as follows:
$$a\cdot_\sigma b=\sigma(a_{(1)}, b_{(1)})a_{(2)}b_{(2)}\sigma^{-1}(a_{(3)}, b_{(3)}),\quad a, b\in H.$$

\begin{prop}\cite{HuPeiCA 2012}
As Hopf algebras, $\mathcal{A}_{r,s}$ can be realized as the cocycle twist $(\mathcal{A}_1)^\sigma$.
The $2$-cocycle $\sigma\colon \mathcal{A}_1\otimes \mathcal{A}_1\to \mathbb{C}$ is given by
\begin{equation*}\label{equ cocycle}
\sigma(x,y)=
\begin{cases}
r^{\frac{\langle \nu,\mu\rangle}{2}} s^{-\frac{\langle \mu,\nu\rangle}{2}}
& \text{if } x=\omega_\mu \text{ or } x=\omega'_\mu,\ y=\omega_\nu \text{ or } y=\omega'_\nu;\\[4pt]
0 & \text{otherwise},
\end{cases}
\end{equation*}
where $\langle\cdot,\cdot\rangle$ denotes the Euler form introduced in \cite{HuPei2008}. The pairings $\langle n-1,n\rangle$ and $\langle n,n-1\rangle$ are to be redefined according to \eqref{equ redefi}.
\end{prop}

\begin{theorem}\label{thm diff factors}
Let \(V_1\) and \(V_2\) be respectively a one-dimensional and an
infinite-dimensional simple \(\mathcal B_{rs}\)-module, and set
$M(V_i)=\mathcal A_{r,s}\otimes_{\mathcal B_{rs}}V_i,\;i=1,2$.
Then $$\operatorname{Hom}_{\mathcal A_{r,s}}(M(V_1),M(V_2))=0.$$
\end{theorem}

\begin{proof}
By part (1) of the proof of Theorem \ref{thm type D isoclass},
$e_{n-1}$ and $e_n$ act invertibly on $V_2$,
while at least one of these two generators acts trivially on $V_1$.
We may assume, without loss of generality, that \(e_nV_1=0\).
Let \(\mathfrak B_{r,s}\) denote the extended Borel algebra,
generated by $e_{n-1}, e_n, \omega_{n-1}, \omega_n, \omega'_{n-1}$ and $\omega'_n$.
Thus \(\mathcal B_{rs}\subseteq \mathfrak B_{r,s}\subseteq \mathcal A_{r,s}\).
Put $N_i=\mathfrak B_{r,s}\otimes_{\mathcal B_{r,s}}V_i,\;i=1,2.$
By transitivity of induction, we have
$M(V_i)=\mathcal A_{r,s}\otimes_{\mathcal B_{r,s}}V_i\simeq \mathcal A_{r,s}\otimes_{\mathfrak B_{r,s}}N_i.$
Therefore Frobenius reciprocity gives
\[
\operatorname{Hom}_{\mathcal A_{r,s}}(M(V_1),M(V_2))
\simeq
\operatorname{Hom}_{\mathfrak B_{r,s}}
\left(N_1, \operatorname{Res}^{\mathcal A_{r,s}}_{\mathfrak B_{r,s}}M(V_2)\right).
\]
It is enough to show that the Hom-space on the right-hand side is zero.

First we observe that \(e_nN_1=0\).
Indeed, since $e_n$ quasicommute with the other generators of $\mathfrak B_{r,s}$,
\(e_n\mathfrak B_{r,s}\subseteq \mathfrak B_{r,s}e_n\).
Thus, for \(b\in\mathfrak B_{r,s}\) and \(v\in V_1\), we have
$e_n(b\otimes v)=b'e_n\otimes v=b'\otimes e_nv=0$
for some \(b'\in\mathfrak B_{r,s}\).

Next we claim that \(e_n\) acts injectively on \(\operatorname{Res}^{\mathcal A_{r,s}}_{\mathfrak B_{r,s}}M(V_2)\).
Since the group-like elements quasicommute with \(e_n\),
the module \(N_2\) is a direct sum of copies of \(V_2\) on which the action of \(e_n\) differs from its action on \(V_2\) only by nonzero scalars.
Since \(e_n\) acts injectively on \(V_2\), it acts injectively on \(N_2\).

Now use the PBW basis of \(\mathcal A_{r,s}\).
Let \(\mathcal P\) be the linear span of the remaining PBW monomials.
Thus \(\mathcal P\) contains the negative root vectors,
the positive rootvectors and group-like elements not contained in \(\mathfrak B_{r,s}\).
Multiplication gives a vector space isomorphism
$\mathcal A_{r,s}\simeq \mathcal P\otimes \mathfrak B_{r,s}.$
Hence $M(V_2) \simeq\mathcal A_{r,s}\otimes_{\mathfrak B_{r,s}}N_2\simeq\mathcal P\otimes N_2$ as vector spaces.

Let \(\{P_{\mathbf a}\}\) be the ordered PBW basis of \(\mathcal P\).
The order is chosen as follows: the negative root vectors are placed first,
then the group-like elements, and finally the positive root vectors. More
precisely, each basis element is ordered in the form
\begin{equation}\label{equ order}
P_a=
\left(\prod_{\alpha\in\Phi^+}^{\succ}\mathcal F_\alpha^{m_\alpha}\right)
\left(\prod_{i=1}^{n-2}\omega_i^{p_i}
      \prod_{i=1}^{n-2}(\omega'_i)^{q_i}\right)
\left(\prod_{\alpha\in\Phi_{\mathcal P}^+}^{\prec}\mathcal E_\alpha^{\ell_\alpha}\right),
\end{equation}
where \(m_\alpha,\ell_\alpha\in\mathbb N\) and \(p_i,q_i\in\mathbb Z\).
Here the orders on the negative and positive parts are opposite to each
other.

Every
nonzero element \(u\in M(V_2)\) can be written uniquely as
$u=\sum_{\mathbf a}P_{\mathbf a}\otimes n_{\mathbf a}$,
where \(n_{\mathbf a}\in N_2\) and only finitely many \(n_{\mathbf a}\) are
nonzero. Choose a maximal index \(\mathbf a_0\) such that
\(n_{\mathbf a_0}\neq 0\).
By the PBW commutation relations, left
multiplication by \(e_n\) has the triangular form (the precise verification can be seen in Lemma \ref{lemm triangular} in Appendix):
\begin{equation*}\label{}
  e_nP_{\mathbf a}=c_{\mathbf a}P_{\mathbf a}e_n+\sum_{\mathbf b<\mathbf a}P_{\mathbf b}b_{\mathbf b},
\end{equation*}
where \(c_{\mathbf a}\in\mathbb C^\times\) and
\(b_{\mathbf b}\in\mathfrak B_{r,s}\).
Therefore the highest PBW term of \(e_nu\) is
$c_{\mathbf a_0}P_{\mathbf a_0}\otimes e_nn_{\mathbf a_0}.$
This term is nonzero because \(c_{\mathbf a_0}\neq 0\) and \(e_n\) acts
injectively on \(N_2\). It cannot be cancelled by lower PBW terms. Hence
\(e_nu\neq 0\), and consequently
\[
\ker\left(
e_n:
\operatorname{Res}^{\mathcal A_{r,s}}_{\mathfrak B_{r,s}}M(V_2)
\to
\operatorname{Res}^{\mathcal A_{r,s}}_{\mathfrak B_{r,s}}M(V_2)
\right)=0.
\]

Let $\varphi\in\operatorname{Hom}_{\mathfrak B_{r,s}}
\left(N_1, \operatorname{Res}^{\mathcal A_{r,s}}_{\mathfrak B_{r,s}}M(V_2)\right)$.
Since \(e_nN_1=0\), for every \(w\in N_1\) we have
$e_n\varphi(w)=\varphi(e_nw)=0$.
Thus \(\varphi(w)\) lies in the kernel of \(e_n\) on
\(\operatorname{Res}^{\mathcal A_{r,s}}_{\mathfrak B_{r,s}}M(V_2)\), which is
zero by the previous paragraph. Hence \(\varphi=0\). Therefore
$\operatorname{Hom}_{\mathcal A_{r,s}}(M(V_1),M(V_2))=0$.
\end{proof}

The above theorem suggests that the induced modules from $\mathcal{B}_{r,s}$ already capture certain intrinsic features of the module category of the whole two-parameter quantum group.
This phenomenon indicates a substantial difference between the representation theories of two-parameter and one-parameter quantum groups.
Thus we conjecture that
\begin{conj}\label{thm Ars}
  $\mathcal{A}_{r,s}$ is not isomorphic to $\mathcal{A}_1$ as $\mathbb C$-algebras.
\end{conj}

\appendix

\section{Additional examples of the type $D$ phenomenon}
\label{app type D examples}

We verify the relations in Table~\ref{tab type D examples}
under the standing assumptions on $r$ and $s$.
The notation for each type follows the corresponding reference.
We recall the root-vector definitions needed below.

We first establish a local identity used in finite
and affine type $A$.

\begin{lemm}\label{lem local A3 commutation}
  Let $a,b,c$ be elements of a $\mathbb{C}$-algebra satisfying
  $ac=ca$ and
  \begin{align}
    ab^2-(r+s)bab+rs\,b^2a&=0,
    \label{eq local A3 serre a}\\
    b^2c-(r+s)bcb+rs\,cb^2&=0.
    \label{eq local A3 serre c}
  \end{align}
  Assume $r+s\neq0$, and set
  \[
    X=ab-rba,\qquad Y=bc-scb.
  \]
  Then
  \[
    XY=rs\,YX.
  \]
\end{lemm}

\begin{proof}
  Using $ac=ca$, we obtain
  \begin{align*}
    D:=XY-rs\,YX
    ={}&ab^2c-sabcb-rbabc\\
       &+r^2s\,bcba+rs^2cbab-r^2s^2cb^2a.
  \end{align*}
  Set
  \[
    S_a=ab^2-(r+s)bab+rs\,b^2a,\qquad
    S_c=b^2c-(r+s)bcb+rs\,cb^2.
  \]
  A direct expansion gives
  \[
    (r+s)D
    =rS_ac+s\,aS_c-rs^2cS_a-r^2sS_ca.
  \]
  Since $S_a=S_c=0$ and $r+s\neq0$, we have $D=0$.
\end{proof}

\medskip
\noindent\textbf{Finite type $A_n$.}
Assume $n\geq3$, and set
\[
  X=e_1e_2-r e_2e_1,\qquad
  Y=e_2e_3-s e_3e_2.
\]
The defining relations give $e_1e_3=e_3e_1$ and
\begin{align*}
  e_1e_2^2-(r+s)e_2e_1e_2+rs\,e_2^2e_1&=0,\\
  e_2^2e_3-(r+s)e_2e_3e_2+rs\,e_3e_2^2&=0.
\end{align*}
Thus the hypotheses of Lemma~\ref{lem local A3 commutation}
hold for $a=e_1$, $b=e_2$, and $c=e_3$.
The lemma therefore gives
\[
  XY=rs\,YX.
\]

In the notation of \cite{BW FIC 2004},
\[
  \mathcal{E}_{2,1}=e_2e_1-r^{-1}e_1e_2,
  \qquad X=-r\mathcal{E}_{2,1}.
\]
Substituting this expression for $X$ and cancelling $-r$,
we obtain
\begin{equation}\label{eq finite A type D example}
  \mathcal{E}_{2,1}Y=rs\,Y\mathcal{E}_{2,1}.
\end{equation}

For the root-vector choice, take the letter order
$e_1<e_3<e_2$ in the local $A_3$ subsystem.
The two-letter good Lyndon words $e_1e_2$ and $e_3e_2$
have quantum brackets
\[
  e_1e_2-r e_2e_1=X,\qquad
  e_3e_2-s^{-1}e_2e_3=-s^{-1}Y.
\]
Hence $\mathcal{E}_{2,1}$ and $Y$ are nonzero scalar multiples
of Lyndon root vectors of degrees
$\alpha_1+\alpha_2$ and $\alpha_2+\alpha_3$, respectively.

\medskip
\noindent\textbf{Finite type $B_n$.}
Assume $n\geq3$, and use the root-vector conventions
of \cite{HuWang2010}, where $E_{i,j'}$ corresponds to
the positive root $\varepsilon_i+\varepsilon_j$.
The recursive definitions give
\begin{align*}
  E_{n-1,n}
  &=e_{n-1}e_n-r^2e_ne_{n-1},\\
  E_{n-2,n}
  &=e_{n-2}E_{n-1,n}-r^2E_{n-1,n}e_{n-2},\\
  E_{n-2,n'}
  &=E_{n-2,n}e_n-rs\,e_nE_{n-2,n},\\
  E_{n-2,(n-1)'}
  &=E_{n-2,n'}e_{n-1}-s^{-2}e_{n-1}E_{n-2,n'}.
\end{align*}
By Lemma~3.2(4) and Lemma~3.4(5) of \cite{HuWang2010},
\[
  E_{n-2,n'}e_n=s^2e_nE_{n-2,n'},\qquad
  E_{n-2,n'}E_{n-1,n}=E_{n-1,n}E_{n-2,n'}.
\]
Using these relations and the recursive definitions, we obtain
\begin{align*}
  &E_{n-2,(n-1)'}e_n-(rs)^2e_nE_{n-2,(n-1)'}\\
  &\quad=
    E_{n-2,n'}(e_{n-1}e_n-r^2e_ne_{n-1})
    -(e_{n-1}e_n-r^2e_ne_{n-1})E_{n-2,n'}\\
  &\quad=
    E_{n-2,n'}E_{n-1,n}-E_{n-1,n}E_{n-2,n'}\\
  &\quad=0.
\end{align*}
Therefore,
\begin{equation}\label{eq finite B type D example}
  E_{n-2,(n-1)'}e_n
  =(rs)^2e_nE_{n-2,(n-1)'}.
\end{equation}

\medskip
\noindent\textbf{Finite type $C_n$.}
Assume $n\geq3$, and use the conventions of \cite{Chen Phd}.
The simple roots are
\[
  \alpha_i=\varepsilon_i-\varepsilon_{i+1}
  \quad (1\leq i<n),\qquad
  \alpha_n=2\varepsilon_n.
\]
Here $E_{i,j'}$ denotes the root vector of degree
$\varepsilon_i+\varepsilon_j$.
The root vector used in the table is determined by
\begin{align*}
  E_{n-2,n-1}
  &=e_{n-2}e_{n-1}-s e_{n-1}e_{n-2},\\
  E_{n-2,n}
  &=E_{n-2,n-1}e_n-s^2e_nE_{n-2,n-1},\\
  E_{n-2,(n-1)'}
  &=E_{n-2,n}e_{n-1}
    -r^{-1}e_{n-1}E_{n-2,n}.
\end{align*}
Thus
\[
  \deg E_{n-2,(n-1)'}
  =\varepsilon_{n-2}+\varepsilon_{n-1},
  \qquad
  \deg e_n=2\varepsilon_n.
\]

Lemma~3.5(4) of \cite{Chen Phd} gives
\[
  E_{i,(n-1)'}e_n=r^2s^2e_nE_{i,(n-1)'},
  \qquad 1\leq i<n-1.
\]
Taking $i=n-2$, we obtain
\begin{equation}\label{eq finite C type D example}
  E_{n-2,(n-1)'}e_n
  =(rs)^2e_nE_{n-2,(n-1)'}.
\end{equation}

\medskip
\noindent\textbf{Finite type $D_n$.}
For $n\geq4$, the simple root generators corresponding to the two
terminal nodes $n-1$ and $n$ attached to $n-2$ satisfy
$e_{n-1}e_n=rs\,e_ne_{n-1}$ in the conventions of \cite{BGH 2006}.

\medskip
\noindent\textbf{Finite types $E_6$, $E_7$, and $E_8$.}
We use the node numbering and root-vector conventions
of \cite{Bai}.
Only the local subdiagram
\[
  2\;-\;4\;-\;3
\]
is needed: nodes $2$ and $3$ are both adjacent to $4$
and are not adjacent to each other.
Set
\[
  a=e_2,\qquad b=e_4,\qquad c=e_3.
\]
The root vector denoted by $[E_{243}]$ in \cite{Bai}
is written here as $\mathcal{E}_{243}$.
Explicitly,
\[
  H:=\mathcal{E}_{243}
  =(ab-sba)c-r^{-1}c(ab-sba),
  \qquad
  \deg H=\alpha_2+\alpha_3+\alpha_4.
\]
The subscript $243$ records the bracketed word
corresponding to this definition.

The defining relations give $ac=ca$ and
\begin{align}
  T_a:=b^2a-(r^{-1}+s^{-1})bab+(rs)^{-1}ab^2&=0,
  \label{eq E local serre a}\\
  T_c:=b^2c-(r^{-1}+s^{-1})bcb+(rs)^{-1}cb^2&=0.
  \label{eq E local serre c}
\end{align}
Expanding $H$ and using $ac=ca$, we obtain
\[
  H=abc-sbac-r^{-1}acb+r^{-1}s\,cba.
\]
Hence
\begin{align*}
  D:=Hb-rs\,bH
  ={}&abcb-r^{-1}acb^2+r^{-1}s\,cbab\\
     &-rs\,babc+rs^2b^2ac-s^2bcba.
\end{align*}
A direct expansion gives
\[
  (r+s)D
  =-rs\,aT_c-s^2cT_a+(rs)^2T_ac+rs^3T_ca.
\]
Since $T_a=T_c=0$ and $r+s\neq0$, we have $D=0$.
Therefore,
\begin{equation}\label{eq finite E type D example}
  \mathcal{E}_{243}e_4=rs\,e_4\mathcal{E}_{243}.
\end{equation}
The calculation applies to $E_6$, $E_7$, and $E_8$,
since it uses only their common subdiagram on $2,3,4$
and the corresponding defining relations.

\medskip
\noindent\textbf{Finite type $F_4$.}
We use the notation of \cite{CHW F4}, where the simple root
generators are denoted by $E_i$.
The subscript in $E_{12332}$ records its defining Lyndon word;
in particular,
\[
  \deg E_{12332}=\alpha_1+2\alpha_2+2\alpha_3.
\]
Lemma~3.3(8) in that reference states
\[
  [E_{12332},E_3]=0,
\]
where the bracket is the quantum bracket.

Write the structure constants as
\[
  \omega_iE_j\omega_i^{-1}=a_{ij}E_j.
\]
For homogeneous root vectors $U,V$ of degrees
$\sum_i u_i\alpha_i$ and $\sum_j v_j\alpha_j$,
the bracket convention is
\[
  [U,V]=UV-p_{U,V}VU,\qquad
  p_{U,V}=\prod_{i,j}a_{ij}^{u_iv_j}.
\]
The structure constant matrix gives
\[
  a_{13}=1,\qquad a_{23}=s^2,\qquad a_{33}=rs^{-1}.
\]
Therefore,
\[
  p_{E_{12332},E_3}
  =a_{13}a_{23}^2a_{33}^2
  =s^4(rs^{-1})^2
  =(rs)^2.
\]
Thus the cited quantum-bracket relation gives
\begin{equation}\label{eq finite F type D example}
  E_{12332}E_3=(rs)^2E_3E_{12332}.
\end{equation}

\medskip
\noindent\textbf{Affine type $A_m^{(1)}$.}
Assume $m\geq3$, and use the conventions of \cite{HRZ 2008}
for $U_{r,s}(\widehat{\mathfrak{sl}}_{m+1})$.
Here $0$ is the affine node, and the nodes $m,0,1$
form the finite-type $A_3$ subdiagram
\[
  m\;-\;0\;-\;1.
\]
Set
\[
  a=e_m,\qquad b=e_0,\qquad c=e_1.
\]
The defining relations give $ac=ca$ and
\[
  ab^2-(r+s)bab+rs\,b^2a=0,\qquad
  b^2c-(r+s)bcb+rs\,cb^2=0.
\]
Therefore Lemma~\ref{lem local A3 commutation} gives
\[
  XY=rs\,YX,
  \qquad
  X=e_me_0-r e_0e_m,\qquad
  Y=e_0e_1-s e_1e_0.
\]

Let
\[
  \beta=\alpha_m+\alpha_0,\qquad
  \eta=\alpha_1+\alpha_0,
\]
and define
\[
  E_{\beta}=e_me_0-r e_0e_m,\qquad
  E_{\eta}=e_1e_0-s^{-1}e_0e_1.
\]
Since $X=E_{\beta}$ and $Y=-sE_{\eta}$,
the preceding identity becomes
\begin{equation}\label{eq affine A type D example}
  E_{\beta}E_{\eta}=rs\,E_{\eta}E_{\beta}.
\end{equation}

Both $\beta$ and $\eta$ are positive real roots.
For the root-vector choice, take the letter order
$e_m<e_1<e_0$ in the local $A_3$ subsystem.
The two-letter good Lyndon words are $e_me_0$ and $e_1e_0$.
By the pairing in \cite[Proposition~2.4]{HRZ 2008},
\[
  \langle\omega_m',\omega_0\rangle=r^{-1},\qquad
  \langle\omega_1',\omega_0\rangle=s.
\]
Their Lyndon brackets, with the convention
of \cite[Definition~3.6]{HRZ 2008}, are therefore
$E_{\beta}$ and $E_{\eta}$.
This choice is made in the local finite-type subsystem,
rather than in the particular affine root-vector system
fixed in that reference.
The restriction $m\geq3$ is needed for this $A_3$ subdiagram.

\medskip
\noindent\textbf{The remaining affine examples.}
We use the affine node numbering in the references
listed in Table~\ref{tab type D examples}.
The elements $e_i$ are the corresponding simple root generators,
with $e_0$ the affine generator.
For the last row of the table, take
\[
  (i,j)=
  \begin{cases}
    (0,1),
      & \text{for } B_n^{(1)},\ D_n^{(1)},\ E_6^{(1)},\\
    (n,0),
      & \text{for } C_n^{(1)},\\
    (0,2),
      & \text{for } F_4^{(1)},\\
    (2,0),
      & \text{for } G_2^{(1)}.
  \end{cases}
\]
Here $n\geq3$ for $B_n^{(1)}$ and $C_n^{(1)}$,
and $n\geq4$ for $D_n^{(1)}$.
In each case,
\[
  a_{ij}=0,\qquad
  \langle\omega_j',\omega_i\rangle=rs.
\]
The pairing coefficient is characterized by
\[
  \omega_i e_j\omega_i^{-1}
  =\langle\omega_j',\omega_i\rangle e_j.
\]
For $D_n^{(1)}$ and $E_6^{(1)}$, these values follow from
the structure constant matrices in \cite{HuZhang2014};
for the other types, see \cite[Chapter~2]{Zhang phd 2007}.

Since $a_{ij}=0$, the corresponding Serre relation gives
\begin{align*}
  0
  &=(\operatorname{ad}_l e_i)(e_j)\\
  &=e_ie_j-\omega_i e_j\omega_i^{-1}e_i\\
  &=e_ie_j-\langle\omega_j',\omega_i\rangle e_je_i\\
  &=e_ie_j-rs\,e_je_i.
\end{align*}
Thus
\begin{equation}\label{eq affine simple type D examples}
  e_ie_j=rs\,e_je_i.
\end{equation}

Under the specialization $r=q$, $s=q^{-1}$,
both $rs$ and $(rs)^2$ become $1$.
All the relations in Table~\ref{tab type D examples}
therefore become commutative relations.
\section{Lyndon root vectors of type $D_n$} \label{appA}

In this appendix we recall the Lyndon root vectors of type $D_n$.

We use the notation of Section \ref{sec pre}.
Let \(U^+_{r,s}(D_n)\) be the positive part of \(U_{r,s}(\mathfrak{so}_{2n})\).
Write \(\mathcal E_{\alpha_i}=e_i\) for \(1\leq i\leq n\).
For each non-simple positive root \(\gamma\in\Phi^+\),
let \((\alpha,\beta)\) be its costandard pair,
i.e. the unique pair of positive roots such that \(\gamma=\alpha+\beta\) and \(\alpha\prec\gamma\prec\beta\), with \(\alpha\) maximal.
If \(\mathcal E_\alpha\in U^+_\zeta\) and \(\mathcal E_\beta\in U^+_\eta\), define
the positive root vector recursively by
$$\mathcal E_\gamma
:=[\mathcal E_\alpha,\mathcal E_\beta]_{\langle\omega'_\eta,\omega_\zeta\rangle}
=\mathcal E_\alpha\mathcal E_\beta-\langle\omega'_\eta,\omega_\zeta\rangle \mathcal E_\beta\mathcal E_\alpha.$$
Equivalently, in type \(D_n\) the positive Lyndon root vectors may be written as follows:
\[
\begin{aligned}
\mathcal E_{i,j} &= e_i\mathcal E_{i+1,j}-s\mathcal E_{i+1,j}e_i,
&& 1\leq i<j\leq n-1,\\
\mathcal E_{n-2,n} &= e_{n-2}e_n-se_ne_{n-2},\\
\mathcal E_{i,n} &= e_i\mathcal E_{i+1,n}-s\mathcal E_{i+1,n}e_i,
&& 1\leq i\leq n-3,\\
\mathcal E_{i,\bar j} &= e_i\mathcal E_{i+1,\bar j}-s\mathcal E_{i+1,\bar j}e_i,
&& 1\leq i<j-1<n-1,\\
\mathcal E_{i,\overline{i+1}} &= \mathcal E_{i,i+2}e_{i+1}
-r^{-1}e_{i+1}\mathcal E_{i,i+2},
&& 1\leq i\leq n-2.
\end{aligned}
\]
For example, in type \(D_4\),
$\mathcal E_{1\overline{3}}, \mathcal E_{1\overline{2}}$, and $\mathcal E_{2\overline{3}}$
denote the positive root vectors corresponding respectively to the positive roots
$\alpha_1+\alpha_2+\alpha_3+\alpha_4, \alpha_1+2\alpha_2+\alpha_3+\alpha_4$, and $\alpha_2+\alpha_3+\alpha_4$.

Similarly, let \(U^-_{r,s}(D_n)\) be the negative part of \(U_{r,s}(\mathfrak{so}_{2n})\).
Write \(\mathcal F_{\alpha_i}=f_i\) for \(1\leq i\leq n\).
For the same costandard pair \((\alpha,\beta)\) of a non-simple positive root \(\gamma=\alpha+\beta\),
the corresponding negative root vector is defined recursively in the opposite order. If
\(\mathcal F_\alpha\in U^-_\zeta\) and \(\mathcal F_\beta\in U^-_\eta\), set
\[
\mathcal F_\gamma
:=[\mathcal F_\beta,\mathcal F_\alpha]_{\langle\omega'_\zeta,\omega_\eta\rangle^{-1}}
=\mathcal F_\beta\mathcal F_\alpha
-\langle\omega'_\zeta,\omega_\eta\rangle^{-1}
\mathcal F_\alpha\mathcal F_\beta .
\]

According to \cite{HRZ 2008}, there is a $\mathbb{Q}$-algebra anti-automorphism $\tau$ of $U_{r,s}(\mathfrak{g})$
such that $\tau(r)=s$, $\tau(s)=r$ and
\begin{equation}\label{equ tau}
  \tau(e_i)=f_i,\qquad
  \tau(f_i)=e_i,\qquad
  \tau(\omega_i)=\omega_i',\qquad
  \tau(\omega_i')=\omega_i,
\end{equation}
we have
$\tau(\langle \omega_i',\omega_j\rangle^{\pm})=\langle \omega_j',\omega_i\rangle^{\mp}$,
and $\tau(\mathcal{E}_\alpha)=\mathcal{F}_\alpha$ for all $\alpha\in\Phi^+$.

\begin{theorem}\cite{Bai}\label{thm convex PBW Lyndon bases}
The sets
\(\{\prod_{\alpha\in\Phi^+}^{\prec}\mathcal E_\alpha^{\,n_\alpha}
\mid n_\alpha\in\mathbb N\}\) and
\(\{\prod_{\alpha\in\Phi^+}^{\succ}\mathcal F_\alpha^{\,n_\alpha}
\mid n_\alpha\in\mathbb N\}\) are convex PBW-type Lyndon bases for
\(U^+_{r,s}(D_n)\) and \(U^-_{r,s}(D_n)\), respectively.
\end{theorem}

For \(n=4\), the Lyndon root vectors are ordered as
\begin{equation}\label{equ convex}
  e_1\prec\mathcal E_{12}\prec\mathcal E_{13}\prec\mathcal E_{14} \prec\mathcal E_{1\bar{3}}\prec\mathcal E_{1\bar{2}}
\prec e_2\prec\mathcal E_{23}\prec\mathcal E_{24}\prec\mathcal E_{2\bar{3}}\prec e_3\prec e_4 .
\end{equation}
\begin{lemm}\cite{Bai}
  $$\{e_1^{a_1}\mathcal{E}_{12}^{a_2}\mathcal{E}_{13}^{a_3}\mathcal{E}_{14}^{a_4}\mathcal{E}_{1\bar{3}}^{a_5}\mathcal{E}_{1\bar{2}}^{a_6}e_2^{a_7}\mathcal{E}_{23}^{a_8}\mathcal{E}_{24}^{a_9}\mathcal{E}_{2\bar{3}}^{a_{10}}e_3^{a_{11}}e_4^{a_{12}}\mid a_i\in \mathbb{Z}^{\geq 0},\, 1\leq i\leq 12\} $$ forms a Lyndon basis of $U^+_{r, s}(D_4)$.
\end{lemm}

\begin{lemm}\label{lem D4}
  In $U^+_{r, s}(D_4)$, we have the following relations:
  \begin{align}
     &e_1\mathcal{E}_{12}=r\mathcal{E}_{12}e_1,  \label{equ D4 1} \\
     & \mathcal{E}_{12}\mathcal{E}_{13}=r\mathcal{E}_{13}\mathcal{E}_{12}, \label{equ D4 2}\\
     & \mathcal{E}_{13}\mathcal{E}_{14}=rs\mathcal{E}_{14}\mathcal{E}_{13}, \label{equ D4 3} \\
     & \mathcal{E}_{14}\mathcal{E}_{1\bar{3}}=s^{-1}\mathcal{E}_{1\bar{3}}\mathcal{E}_{14} \label{equ D4 4}\\
     &  \mathcal{E}_{1\bar{3}}\mathcal{E}_{1\bar{2}}=s^{-1}\mathcal{E}_{1\bar{2}}\mathcal{E}_{1\bar{3}}, \label{equ D4 5}\\
     & \mathcal{E}_{1\bar{2}}e_2=s^{-1}e_2\mathcal{E}_{1\bar{2}},  \label{equ D4 6}\\
     & e_2\mathcal{E}_{23}=r\mathcal{E}_{23}e_2, \label{equ D4 7}  \\
     & \mathcal{E}_{23}\mathcal{E}_{24}=rs\mathcal{E}_{24}\mathcal{E}_{23}, \label{equ D4 8}\\
     & \mathcal{E}_{24}\mathcal{E}_{2\bar{3}}=s^{-1}\mathcal{E}_{2\bar{3}}\mathcal{E}_{24},  \label{equ D4 9} \\
     & \mathcal{E}_{2\bar{3}}e_3=s^{-1}e_3\mathcal{E}_{2\bar{3}}, \label{equ D4 10}\\
     & e_3e_4=rse_4e_3. \label{equ D4 11}
  \end{align}
\end{lemm}

\begin{proof}
 (\ref{equ D4 1}) and (\ref{equ D4 7}) are directly from the Serre relations.
 (\ref{equ D4 6}) and (\ref{equ D4 10}) are due to Lemma 3.16 (4) in \cite{Bai}.
 (\ref{equ D4 11}) is from the definition of $U_{r, s}(\mathfrak{so}_{2n})$.
 The remaining identities follow from Equation (101) of Theorem $D_n$ in \cite{Kha Pac02}.
\end{proof}

The following calculations will be useful in Appendix \ref{appB}.
Applying $\phi$ to (\ref{equ D4 4}), we obtain that $\mathcal{E}_{13}\omega_3^{-1}\omega_4\mathcal{E}_{1\bar{3}}=s^{-1}\mathcal{E}_{1\bar{3}}\mathcal{E}_{13}\omega_3^{-1}\omega_4$.
That is,
\begin{equation}\label{equ 13 13bar}
  \mathcal{E}_{13}\mathcal{E}_{1\bar{3}}=r^2s\mathcal{E}_{1\bar{3}}\mathcal{E}_{13}.
\end{equation}
Similarly, we can verify that
\begin{equation}\label{equ 12 14}
  \mathcal{E}_{12}\mathcal{E}_{14}=r\mathcal{E}_{14}\mathcal{E}_{12},\qquad e_2\mathcal{E}_{24}=r\mathcal{E}_{24}e_2,
\end{equation}
\begin{equation}\label{equ 23 23bar}
  \mathcal{E}_{23}\mathcal{E}_{2\bar{3}}=r^2s\mathcal{E}_{2\bar{3}}\mathcal{E}_{23}, \qquad \mathcal{E}_{2\bar{3}}e_4=r^2se_4\mathcal{E}_{2\bar{3}}.
\end{equation}

We calculate that
\begin{align*}
  \mathcal{E}_{13}\mathcal{E}_{24} & =\, e_1\mathcal{E}_{23}\mathcal{E}_{24}-s\mathcal{E}_{23}e_1\mathcal{E}_{24} \\
   & =\,rse_1\mathcal{E}_{24}\mathcal{E}_{23}-s^2\mathcal{E}_{23}\mathcal{E}_{24}e_1-s\mathcal{E}_{23}\mathcal{E}_{14}  \\
   & =\,rs\mathcal{E}_{14}\mathcal{E}_{23}+rs^2\mathcal{E}_{24}e_1\mathcal{E}_{23}-s^2\mathcal{E}_{23}\mathcal{E}_{24}e_1-s\mathcal{E}_{23}\mathcal{E}_{14} \\
   & =\,rs^2\mathcal{E}_{24}e_1\mathcal{E}_{23}-rs^3\mathcal{E}_{24}\mathcal{E}_{23}e_1+rs\mathcal{E}_{14}\mathcal{E}_{23}-s\mathcal{E}_{23}\mathcal{E}_{14} \\
   & =\,rs^2\mathcal{E}_{24}\mathcal{E}_{13}+rs\mathcal{E}_{14}\mathcal{E}_{23}-s\mathcal{E}_{23}\mathcal{E}_{14},
\end{align*}
and
\begin{align*}
  \mathcal{E}_{14}\mathcal{E}_{23} & =\, \mathcal{E}_{14}(e_2e_3-se_3e_2) \\
   &=\,e_2\mathcal{E}_{14}e_3-s\mathcal{E}_{1\bar{3}}e_2-r^{-1}se_3\mathcal{E}_{14}e_2  \\
   &=\, e_2\mathcal{E}_{1\bar{3}}+r^{-1}e_2e_3\mathcal{E}_{14}-s\mathcal{E}_{1\bar{3}}e_2-r^{-1}se_3e_2\mathcal{E}_{14}  \\
   &=\,e_2\mathcal{E}_{1\bar{3}}-s\mathcal{E}_{1\bar{3}}e_2+r^{-1}\mathcal{E}_{23}\mathcal{E}_{14}.
\end{align*}
Combining these two equations, we obtain that
\begin{align}
  \mathcal{E}_{13}\mathcal{E}_{24} & =\,rs^2\mathcal{E}_{24}\mathcal{E}_{13}+rse_2\mathcal{E}_{1\bar{3}}-rs^2\mathcal{E}_{1\bar{3}}e_2 \\
   & =\,rs^2\mathcal{E}_{24}\mathcal{E}_{13}+(r^2s-rs^2)\mathcal{E}_{1\bar{3}}e_2-r^2s\mathcal{E}_{1\bar{2}}. \nonumber
\end{align}
We have also calculated that
\begin{equation}\label{}
  \mathcal{E}_{23}e_4=rs^2e_4\mathcal{E}_{23}+rs\mathcal{E}_{2\bar{3}},
\end{equation}
in Proposition \ref{aut D4}, and moreover
\begin{align}\label{equ 13 e4}
  \mathcal{E}_{13}e_4& =\, e_1\mathcal{E}_{23}e_4-s\mathcal{E}_{23}e_1e_4 \\
   & =\,rs^2e_1e_4\mathcal{E}_{23}+rse_1\mathcal{E}_{2\bar{3}}-rs^3e_4\mathcal{E}_{23}e_1-rs^2\mathcal{E}_{2\bar{3}}e_1  \nonumber \\
   & =\,rs^2e_4\mathcal{E}_{13}+rse_1\mathcal{E}_{2\bar{3}}-rs^2\mathcal{E}_{2\bar{3}}e_1  \nonumber \\
   & =\,rs^2e_4\mathcal{E}_{13}+rs\mathcal{E}_{1\bar{3}}. \nonumber
\end{align}

\section{Triangular commutation relations in type $D_4$}\label{appB}

We calculate the triangular commutation relations  used in Section~\ref{sec induced}.
The following elementary formula will be useful in Lemma \ref{lem:e4-negative-computations}.
\begin{lemm} \label{lem:q-commutator-power}
Let \(R\) be an associative algebra, and let \(x,Y,A\in R\).
Suppose that \(q\in\mathbb C\) and
$$xY=Yx+A,\qquad AY=qYA .$$
Then, for every \(m\geq1\), one has
$$xY^m=Y^m x+[m]_qY^{m-1}A,$$
where $[m]_q=1+q+\cdots+q^{m-1}$.
\end{lemm}

\begin{proof}
We prove the formula by induction on \(m\). The case \(m=1\) is exactly the
assumption \(xY=Yx+A\). Suppose that the formula holds for some \(m\geq1\).
Then
\[
\begin{aligned}
xY^{m+1}
&=(Y^m x+[m]_qY^{m-1}A)Y\\
&=Y^m(Yx+A)+[m]_qY^{m-1}AY\\
&=Y^{m+1}x+Y^mA+q[m]_qY^mA\\
&=Y^{m+1}x+(1+q[m]_q)Y^mA\\
&=Y^{m+1}x+[m+1]_qY^mA.
\end{aligned}
\]
Thus the statement follows for all \(m\geq1\).
\end{proof}

\begin{lemm}\label{lem:e4-negative-computations}
For \(m\geq 1\), set $[m]_q:=\frac{1-q^m}{1-q}\;(q\neq 1)$.
Then
\begin{align}
\label{eq:e4-f4-power}
e_4f_4^m
&=f_4^m e_4
+[m]_{rs^{-1}}f_4^{m-1}
\frac{(r^{-1}s)^{m-1}\omega_4-\omega'_4}{r-s},\\
\label{eq:e4-F24-power}
e_4\mathcal F_{24}^m
&=\mathcal F_{24}^m e_4
+[m]_{rs^{-1}}\mathcal F_{24}^{m-1}f_2\omega'_4,\\
\label{eq:e4-F14-power}
e_4\mathcal F_{14}^m
&=\mathcal F_{14}^m e_4
+[m]_{rs^{-1}}\mathcal F_{14}^{m-1}\mathcal F_{12}\omega'_4,\\
\label{eq:e4-F2bar3-power}
e_4\mathcal F_{2\bar 3}^m
&=\mathcal F_{2\bar 3}^m e_4
+[m]_{rs^{-1}}\mathcal F_{2\bar 3}^{m-1}
  \mathcal F_{23}\omega'_4,\\
\label{eq:e4-F1bar3-power}
e_4\mathcal F_{1\bar 3}^m
&=\mathcal F_{1\bar 3}^m e_4
+[m]_{rs^{-1}}\mathcal F_{1\bar 3}^{m-1}
  \mathcal F_{13}\omega'_4,\\
\label{eq:e4-F1bar2-power}
e_4\mathcal F_{1\bar 2}^m
&=\mathcal F_{1\bar 2}^m e_4.
\end{align}
\end{lemm}

\begin{proof}
The first identity is originally from  Lemma 2.6 in \cite{BGH 2007}.

Let \(H_4=(\omega_4-\omega'_4)/(r-s)\), then $e_4f_4=f_4e_4+H_4$.
Since  $\omega_4f_2=rf_2\omega_4$ and $\omega'_4f_2=sf_2\omega'_4$,
we get
\[
\begin{aligned}
e_4\mathcal F_{24}
&=e_4(f_4f_2-rf_2f_4)\\
&=(f_4e_4+H_4)f_2-rf_2(f_4e_4+H_4)\\
&=\mathcal F_{24}e_4+H_4f_2-rf_2H_4\\
&=\mathcal F_{24}e_4+f_2\omega'_4.
\end{aligned}
\]

We now derive the formula for $e_4\mathcal F_{24}^m$.
First note that $f_2\omega'_4\mathcal F_{24}=rs^{-1}\mathcal F_{24}f_2\omega'_4$.
Indeed, \(\mathcal F_{24}f_2=sf_2\mathcal F_{24}\),and
\(\omega'_4\mathcal F_{24}=r\mathcal F_{24}\omega'_4\).
Therefore $f_2\omega'_4\mathcal F_{24}=r f_2\mathcal F_{24}\omega'_4=rs^{-1}\mathcal F_{24}f_2\omega'_4$.

Applying Lemma \ref{lem:q-commutator-power} with \(x=e_4\), \(Y=\mathcal F_{24}\),
\(A=f_2\omega'_4\), and \(q=rs^{-1}\), we obtain
\[
e_4\mathcal F_{24}^m
=\mathcal F_{24}^m e_4
+[m]_{rs^{-1}}\mathcal F_{24}^{m-1}f_2\omega'_4.
\]

The other identities follow from the recursive definitions. For example,
\[
\begin{aligned}
e_4\mathcal F_{14}
&=e_4(\mathcal F_{24}f_1-rf_1\mathcal F_{24})\\
&=(\mathcal F_{24}e_4+f_2\omega'_4)f_1
-rf_1(\mathcal F_{24}e_4+f_2\omega'_4)\\
&=\mathcal F_{14}e_4
+(f_2f_1-rf_1f_2)\omega'_4\\
&=\mathcal F_{14}e_4+\mathcal F_{12}\omega'_4.
\end{aligned}
\]
Similarly, using \(e_4f_3=f_3e_4\) and
\(\omega'_4f_3=rs\,f_3\omega'_4\), we obtain
\[
\begin{aligned}
e_4\mathcal F_{2\bar 3}
&=e_4(f_3\mathcal F_{24}-s^{-1}\mathcal F_{24}f_3)\\
&=f_3(\mathcal F_{24}e_4+f_2\omega'_4)
-s^{-1}(\mathcal F_{24}e_4+f_2\omega'_4)f_3\\
&=(f_3\mathcal F_{24}-s^{-1}\mathcal F_{24}f_3)e_4+\bigl(f_3f_2-rf_2f_3\bigr)\omega'_4\\
&=\mathcal F_{2\bar 3}e_4+\mathcal F_{23}\omega'_4.
\end{aligned}
\]
Using this relation, together with \(e_4f_1=f_1e_4\) and
\(\omega'_4f_1=f_1\omega'_4\), we get
\[
\begin{aligned}
e_4\mathcal F_{1\bar 3}
&=e_4(\mathcal F_{2\bar 3}f_1-rf_1\mathcal F_{2\bar 3})\\
&=(\mathcal F_{2\bar 3}e_4+\mathcal F_{23}\omega'_4)f_1-rf_1(\mathcal F_{2\bar 3}e_4+\mathcal F_{23}\omega'_4)\\
&=(\mathcal F_{2\bar 3}f_1-rf_1\mathcal F_{2\bar 3})e_4+\bigl(\mathcal F_{23}f_1-rf_1\mathcal F_{23}\bigr)\omega'_4\\
&=\mathcal F_{1\bar 3}e_4+\mathcal F_{13}\omega'_4.
\end{aligned}
\]

Moreover, for
\(Y\in\{\mathcal F_{24},\mathcal F_{14},
\mathcal F_{2\bar 3},\mathcal F_{1\bar 3},
\mathcal F_{1\bar 2}\}\), define \(T_Y\) by
\[
T_Y=
\begin{cases}
f_2\omega'_4, & Y=\mathcal F_{24},\\
\mathcal F_{12}\omega'_4, & Y=\mathcal F_{14},\\
\mathcal F_{23}\omega'_4, & Y=\mathcal F_{2\bar 3},\\
\mathcal F_{13}\omega'_4, & Y=\mathcal F_{1\bar 3},\\
0, & Y=\mathcal F_{1\bar 2}.
\end{cases}
\]
Then \(e_4Y=Ye_4+T_Y\).
We next check the commutation relations between \(T_Y\) and \(Y\).

Applying the anti-automorphism $\tau$ (\ref{equ tau}) on (\ref{equ 13 13bar})-(\ref{equ 23 23bar}) gives
\[
\begin{aligned}
&f_2\mathcal F_{24} =s^{-1}\mathcal F_{24}f_2, \hspace{4.8em}
\mathcal F_{12}\mathcal F_{14} =s^{-1}\mathcal F_{14}\mathcal F_{12},\\
&\mathcal F_{23}\mathcal F_{2\bar 3} =r^{-1}s^{-2}\mathcal F_{2\bar 3}\mathcal F_{23}, \qquad
\mathcal F_{13}\mathcal F_{1\bar 3}=r^{-1}s^{-2}\mathcal F_{1\bar 3}\mathcal F_{13}.
\end{aligned}
\]
On the other hand, the toral relations imply
$\omega'_4\mathcal F_{24}=r\mathcal F_{24}\omega'_4,\;
\omega'_4\mathcal F_{14}=r\mathcal F_{14}\omega'_4,\;
\omega'_4\mathcal F_{2\bar 3}=r^2s\,\mathcal F_{2\bar 3}\omega'_4,\;
\omega'_4\mathcal F_{1\bar 3}=r^2s\,\mathcal F_{1\bar 3}\omega'_4.$
Therefore,
$T_YY=rs^{-1}YT_Y$
for
\(Y\in\{\mathcal F_{24},\mathcal F_{14},
\mathcal F_{2\bar 3},\mathcal F_{1\bar 3}\}\). For
\(Y=\mathcal F_{1\bar 2}\), the same relation is trivial because
\(T_{\mathcal F_{1\bar 2}}=0\). Applying Lemma
\ref{lem:q-commutator-power} with \(x=e_4\), \(A=T_Y\), and
\(q=rs^{-1}\), we obtain, for every \(m\geq 1\),
\[
e_4Y^m=Y^m e_4+[m]_{rs^{-1}}Y^{m-1}T_Y.
\]
Then we can get the desired equation.
\end{proof}

\begin{lemm} \label{lemm triangular}
Let \(\mathfrak B_{r,s}\) denote the extended Borel subalgebra of $U_{r, s}(\mathfrak{so}_8)$,
generated by $e_3, e_4, \omega_3, \omega_4, \omega'_3$ and $\omega'_4$.
Let \(\mathcal P\) be the span of the remaining PBW monomials.
Then multiplication gives a vector-space decomposition
$U_{r, s}(\mathfrak{so}_8)\simeq \mathcal P\otimes\mathfrak B_{r,s}$.
If \(P_{\mathbf a}\) is a PBW monomial in \(\mathcal P\), then
\begin{equation}\label{equ triangular}
e_4P_{\mathbf a}=q_{\mathbf a}P_{\mathbf a}e_4+\sum_{\mathbf b\prec\mathbf a}P_{\mathbf b}b_{\mathbf b},
\end{equation}
where \(q_{\mathbf a}\in\mathbb K^\times\), \(b_{\mathbf b}\in\mathfrak B_{r,s}\),
and \(\prec\) is the PBW order on the \(\mathcal P\)-monomials (\ref{equ order}).
\end{lemm}

\begin{proof}
The vector space decomposition follows from the triangular decomposition and
from the PBW bases fixed above.
For negative PBW factors,
we only need to calculate the commutation relations between $e_4$ and those negative composite root vectors containing $f_4$,
these relations are listed in Lemma \ref{lem:e4-negative-computations}.
For group-like elements, the defining relations give only nonzero scalar factors.
For positive PBW factors not belonging to \(\mathfrak B_{r,s}\), we use the fact that (\ref{equ convex}) is indeed convex order.
In each case the leading term is obtained by moving \(e_4\) to the right,
while every extra term has strictly smaller \(\mathcal P\)-part.
This proves the triangular formula (\ref{equ triangular}).
\end{proof}

\bibliographystyle{amsalpha}

\end{document}